\documentclass[12pt,epsf,amssymb]{article}

\usepackage{graphicx}
\usepackage{amsfonts}
\usepackage{amsmath}
\usepackage[usenames]{color}
\usepackage{amssymb}
\usepackage{amsthm}
\usepackage[mathscr]{eucal} 
\usepackage{url}

\def\N{\mathbb{N}}
\def\Z{\mathbb{Z}}
\def\Q{\mathbb{Q}}
\def\R{\mathbb{R}}
\def\C{\mathbb{C}}

\usepackage[%height=22.5cm,
width=16.6cm]{geometry}

\begin{document}

\baselineskip 0.6cm
\newcommand{\nequiv}{\mbox{\ooalign{\hfil/\hfil\crcr$\equiv$}}}
\newcommand{\nsupset}{\mbox{\ooalign{\hfil/\hfil\crcr$\supset$}}}
\newcommand{\nni}{\mbox{\ooalign{\hfil/\hfil\crcr$\ni$}}}
\newcommand{\nin}{\mbox{\ooalign{\hfil/\hfil\crcr$\in$}}}

\newcommand{\vev}[1]{ \left\langle {#1} \right\rangle }
\newcommand{\bra}[1]{ \langle {#1} | }
\newcommand{\ket}[1]{ | {#1} \rangle }
\newcommand{\Dsl}{\mbox{\ooalign{\hfil/\hfil\crcr$D$}}}
\newcommand{\Slash}[1]{{\ooalign{\hfil/\hfil\crcr$#1$}}}
\newcommand{\EV}{ {\rm eV} }
\newcommand{\KEV}{ {\rm keV} }
\newcommand{\MEV}{ {\rm MeV} }
\newcommand{\GEV}{ {\rm GeV} }
\newcommand{\TEV}{ {\rm TeV} }

\def\diag{\mathop{\rm diag}\nolimits}
\def\tr{\mathop{\rm tr}}

\def\Spin{\mathop{\rm Spin}}
\def\SO{\mathop{\rm SO}}
\def\O{\mathop{\rm O}}
\def\SU{\mathop{\rm SU}}
\def\U{\mathop{\rm U}}
\def\Sp{\mathop{\rm Sp}}
\def\SL{\mathop{\rm SL}}

\def\change#1#2{{\color{blue}#1}{\color{red} [#2]}\color{black}\hbox{}}

\theoremstyle{definition}
\newtheorem{thm}{Theorem}[section]
\newtheorem{defn}[thm]{Definition}
\newtheorem{exmpl}[thm]{Example}
\newtheorem{props}[thm]{Proposition}
\newtheorem{lemma}[thm]{Lemma}
\newtheorem{rmk}[thm]{Remark}
\newtheorem{notn}[thm]{Notation}
\newtheorem{anythng}[thm]{}

\begin{titlepage}

%\begin{flushright}
%Version \today
%\end{flushright}

% \begin{flushright}
% IPMUxx-xxxx
% \end{flushright}

\vskip 1cm
\begin{center}
  
{\large \bf A note on varieties of weak CM-type}
 
\vskip 1.2cm

\renewcommand{\thefootnote}{\fnsymbol{footnote}}
Masaki Okada$^a$ and Taizan Watari$^b$
\setcounter{footnote}{0}
\renewcommand{\thefootnote}{\arabic{footnote}}

\vskip 0.4cm

$^a$masaki.okada@ipmu.jp, \quad $^b$taizan.watari@ipmu.jp

\vskip 0.4cm
  
{\it 
     Kavli Institute for the Physics and Mathematics of the Universe (WPI),\\
     University of Tokyo, Kashiwa-no-ha 5-1-5, 277-8583, Japan
}
 
\vskip 1.5cm
    
\abstract{
CM-type projective varieties $X$ of complex dimension $n$ are characterized by their CM-type rational Hodge structures on the cohomology groups.
One may impose such a condition in a weakest form when the canonical bundle of $X$ is trivial;
the rational Hodge structure on the level-$n$ subspace of $H^n(X;\Q)$ is required to be of CM-type.
This brief note addresses the question whether this weak condition implies that the Hodge structure on the entire $H^\ast(X;\Q)$ is of CM-type.
We study in particular abelian varieties when the dimension of the level-$n$ subspace is two or four, and K3$\times T^2$.
It turns out that the answer is affirmative.
Moreover, such an abelian variety is always isogenous to a product of CM-type elliptic curves or abelian surfaces.
This extends a result of Shioda and Mitani in 1974.
} 

\vskip 0.4cm



\end{center}
\end{titlepage}

% \tableofcontents

%%%%%%%%%%%%%%%%%%%%%%%%%%%%%%%%%%%%%%%%%%%%%%%%%
\section{Introduction}
%%%%%%%%%%%%%%%%%%%%%%%%%%%%%%%%%%%%%%%%%%%%%%%%%

The notion of complex multiplication of an elliptic curve $E$
has been generalized to the notion of CM-type of a rational 
Hodge structure on $H^1(E;\Q)$. It is conventional to define whether 
a complex projective variety $X$ is of CM-type or not by whether the rational Hodge 
structure of the cohomology groups of $X$ is of CM-type or not.
There can still be variations, however, on which part of the cohomology group 
$H^*(X;\Q)$ we demand that the rational Hodge structure be of CM-type. 
Typical ones are as follows:
\begin{defn}
\label{defn:strong-CM}
A complex projective non-singular variety $X$ of complex dimension $n$ 
is of {\it strong CM-type} when the polarized rational Hodge structure 
on $H^m(X;\Q)$ is of CM-type for all $m=0,1,\cdots, 2n$
 (see Def. \ref{def:End-CM-Hdg-str} for a rational Hodge structure 
that is of CM-type).
\end{defn}
\begin{defn}
\label{defn:prim-CM}
A complex projective non-singular variety $X$ of complex dimension $n$
is of {\it CM-type} when the polarized rational Hodge structure 
on the primitive subspace $H^n_{\rm prim}(X;\Q)$ of $H^n(X;\Q)$ is of CM-type.
Here, the {\it primitive subspace} is the orthogonal complement of the 
subspace that is generated by the elements in the 
form of $J^k \cdot H^{n-2k}(X;\Q)$ for the polarization $J$ of the projective 
variety $X$ and some integer $0 < k \leq n/2$. 
\end{defn}
We follow \cite[Def.\ 4]{MR3319058} (cf.\ also \cite{rohde2009cyclic}),
make a clear distinction between the Def.\ \ref{defn:strong-CM} from others,
and look into the possible differences between them. 
% MR3319058 = Tretkoff 2014 transcendence and 
% rohde2009 cyclic cover
If a complex projective non-singular $n$-dimensional variety $X$ is of 
strong CM-type, then obviously it is of CM-type. 

In this brief note, the authors are concerned about one more definition:
\begin{defn}
\label{defn:weak-CM}
Let $X$ be a complex projective non-singular variety of complex dimension $n$ 
whose canonical bundle ${\rm det}(T^*X)$ is trivial.\footnote{
The level-$n$ subspace of $H^n(X;\Q)$ is well-defined for a general complex projective non-singular variety $X$ \cite[Def.\ 3.1]{MR3663605},
and the notion of weak CM-type may be extended to such varieties.
We will not explore such possibilities in this note, however.
} % 
Let $[H^n(X;\Q)]_{\ell=n}$ denote the minimum Hodge substructure in $H^n(X;\Q)$ such that $[H^n(X;\Q)]_{\ell=n}\otimes\C$ contains $\Omega^{(n,0)}$, where $\Omega^{(n,0)}$ is a generator of the 
1-dimensional $H^{n,0}(X;\C)$.
We call $[H^n(X;\Q)]_{\ell=n}$ the \textit{level-$n$ subspace}.
The orthogonal component of the level-$n$ subspace within $H^n(X;\Q)$ is denoted by $[H^n(X;\Q)]_{\ell<n}$.

Now, $X$ is said to be of {\it weak CM-type} when the rational 
Hodge structure on the level-$n$ subspace is of CM-type. 
\end{defn}
For a complex projective non-singular $n$-dimensional variety $X$
whose canonical bundle is trivial, $X$ is of weak CM-type if $X$ is of 
strong CM-type (or just CM-type). Led by the question whether the converse 
of this statement holds true (see a note at the end of Introduction for 
a thought from string theory), the authors work on the easiest classes of 
such varieties in this note. Main results are the followings: 
\begin{thm}
\label{thm:intro-main}
{\it Let $A$ be an abelian variety of complex dimension $n$ that is 
of weak CM-type. 
It is then also of strong CM-type, at least if the level-$n$ subspace is 
of dimension $2n' = 2$ or $4$. }
\end{thm}
\begin{thm}
  \label{thm:intro-explct}
{\it Let $A$ be as in Thm.\ \ref{thm:intro-main}. Then $A$ is isogenous 
to the product of $n$ copies of an identical CM elliptic curve $E$ when $n'=1$.
When $n'=2$, $A$ is isogenous to either one of the followings:
\begin{itemize}
\item $(E_1)^r \times (E_2)^{n-r}$ where $E_1$ and $E_2$ are CM elliptic 
curves that are not mutually isogenous, and $0<r<n$, 
\item $S^r$ where $S$ is a simple CM-type abelian surface and $n=2r$.  
\end{itemize}
}
\end{thm}
\begin{thm}
\label{thm:K3xT2}
{\it Let $S$ and $E$ be a complex projective K3 surface and an elliptic curve, 
respectively. If $X = S \times E$ is of weak CM-type, then it is 
also of strong CM-type. }
\end{thm}
\begin{thm}
\label{thm:dcmps}
{\it Let $X_1$ and $X_2$ be complex projective non-singular varieties with 
trivial canonical bundles. When $X_1 \times X_2$ is of weak CM-type, 
then both $X_1$ and $X_2$ are also of weak CM-type. } 
\end{thm}
Proofs of Thms.\ \ref{thm:intro-main} and \ref{thm:intro-explct} are 
found in section \ref{sec:abel}, while Thm.\ \ref{thm:K3xT2} is proved 
in section \ref{sec:K3xT2}; see Rmk. \ref{rmk:dcmps} for Thm. \ref{thm:dcmps}.
 Those proofs only involve elementary 
linear algebra and Galois theory. Thms. \ref{thm:intro-main} and \ref{thm:intro-explct} generalize a part of the results of \cite{shioda1974singular}, 
which dealt with abelian surfaces ($n=2$) with $2n'=2$. 

Section \ref{sec:open-prblm} is expository in nature.
The question whether a weak CM-type abelian variety is of strong
CM-type or not is formulated by using the language of
Dodson \cite{MR735406} in section \ref{ssec:open-abel}, in a way
applicable to the cases with the reflex degree $2n'$ larger than 4. 
A few known examples of Calabi--Yau threefolds of weak CM-type are collected 
in section \ref{ssec:exmpl-CY}. The present authors are physicists
by training; section \ref{sec:open-prblm} is intended to draw
attention of professional mathematicians to this problem\footnote{
Although this still appears to be an open problem in the eyes
of the present authors, it is conceivable that much more
is known and understood than in the survey in this section \ref{sec:open-prblm}.
The present authors would appreciate it very much if they were
introduced to such existing results and observations that are missing
in this note. 
} %
for further progress.
% The authors are happy to follow the advice from the editor and/or referees
% whether to drop section \ref{sec:open-prblm} entirely or to retain it 
% in some form. 

\vspace{3mm}

{\bf A note (some thoughts from string theory):} Here, we leave a brief 
note on a thought in string theory that motivates the converse problem 
we study in this article. The rest of this article is mathematics, 
and is readable without this motivation, however.  

Superstring theory assigns to a Ricci-flat K\"{a}hler manifold $(M, J, ds^2)$ 
with $\dim_\C M \leq 4$ a modular invariant $N=(1,1)$ superconformal 
field theory on Riemann surfaces. The moduli space of $(M, J, ds^2)$ 
is therefore mapped to the moduli space of superconformal field theories 
(SCFTs). Within the latter moduli space are rational SCFTs, where 
the vertex operator superalgebra of an SCFT is maximally enhanced to the 
extent that there is only a finite number of irreducible representations. 
Given the map between the two moduli spaces, one may wonder where 
in the moduli space of $(M,J, ds^2)$ the corresponding SCFTs are rational. 
Study in string theory \cite{Moore:1998pn}, \cite{Moore:1998zu}
in the case $M=T^2$ revealed that the corresponding SCFT is rational 
if and only if the elliptic curve $(M, J)$ is of CM-type, and the mirror 
elliptic curve (determined from the data $(M,J,ds^2)$) is also of CM-type. 
Study of \cite[Thm. 4.5.5]{Wendland:2000ye}, \cite[Thm. 2.5]{Chen:2005gm} 
and \cite[\S3.1]{Kidambi:2022wvh} also hints that we may not 
lose an SCFT that is rational by restricting our attention to 
$(M, J, ds^2)$'s that are projective. 

Ref. \cite{Gukov:2002nw} hints at the possibility that the SCFTs that are 
rational may correspond to $(M,J,ds^2)$ where both $(M,J)$ and the mirror 
manifold $(W,J^\circ)$ of $(M,J,ds^2)$ are {\it something characterized 
by CM-type}, and the endomorphism algebras of $(M,J)$ and $(W,J^\circ)$ 
are {\it somehow related}. It is an on-going research problem in string 
theory to identify the right 
statement (as well as proving it) by sharpening such phrases {\it something 
characterized by} or {\it somehow related}. As studies in the cases with 
higher dimensional $M$ proceed (as in \cite{Chen:2005gm}, 
\cite{Kidambi:2022wvh} and \cite{Okada:2022jnq}), it became clear that 
we need to distinguish multiple different versions of definitions of 
CM-type on varieties (rational Hodge structures) and understand which version 
of the definitions implies other versions and vice versa (beyond what is 
reviewed in the appendix).  

If there is a way to characterize all the data $(M,J,ds^2)$ 
for rational SCFTs all at once, covering all the projective varieties $(M,J)$ 
with a trivial canonical bundle of various topological types, then the 
characterization condition must be about the Hodge structure of some 
universal nature. 
Any projective variety $(M,J)$ of $\dim_\C M=n$ with a trivial canonical 
bundle and its mirror manifold $W$ have their unique level-$n$ Hodge 
substructures, $[H^n(M;\Q)]_{\ell =n}$ in $H^n(M;\Q)$, 
and $[H^n(W;\Q)]_{\ell = n}$ within $H^n(W;\Q)$, and they are both 
simple rational Hodge structures.  
If a universal characterization ever exists, therefore, it is likely 
that the conditions will be in the form that both $(M,J)$ and its mirror 
$(W,J^\circ)$ 
are of weak CM-type, and that there is a Hodge isomorphism between 
the simple rational Hodge structures $[H^n(M;\Q)]_{\ell =n}$ and 
$[H^n(W;\Q)]_{\ell =n}$ (plus one more condition in \cite{Okada:2022jnq}
that is not relevant in this article).

It is also known in the case $(M,J)$ is an abelian variety that 
both $(M,J)$ and $(W,J^\circ)$ are of strong CM-type when the corresponding 
SCFT is rational \cite{Chen:2005gm}. 
It is therefore necessary to dig into the question whether a projective 
variety with a trivial canonical bundle of weak CM-type is also automatically 
of strong CM-type or not. If the answer to this question is negative, 
then the likely universal characterization written above alone may be too 
naive. 

The authors---with background in string theory---will be delighted 
if this note motivates mathematicians to make progress in this converse 
problem.  

%%%%%%%%%%%%%%%%%%%%%%%%%%%%%%%%%%%%%%%%%%%%%%%%%%%%
% \subsection*{Acknowledgements}  % in BrE
\subsection*{Acknowledgments}   % in AmE
%%%%%%%%%%%%%%%%%%%%%%%%%%%%%%%%%%%%%%%%%%%%%%%%%%%%%

The authors are grateful to Yasuhiro Goto and Keita Kanno 
for useful comments and fruitful discussion. 
This work was supported in part by 
FoPM, WINGS Program of the University of Tokyo, 
JSPS Research Fellowship for Young Scientists,
JSPS KAKENHI Grant Number JP23KJ0650 (MO), 
and a Grant-in-Aid for Scientific Research on 
Innovative Areas 6003 and the WPI program (MO and TW), MEXT, Japan.
Declarations of interest: none.

%%%%%%%%%%%%%%%%%%%%%%%%%%%%%%%%%%%%%%%%%%%%%%%%
\section{Abelian Varieties of Weak CM-type}
\label{sec:abel}
%%%%%%%%%%%%%%%%%%%%%%%%%%%%%%%%%%%%%%%%%%%%%%%%%%

Let $A$ be an abelian variety of dimension $n$.
Let its rational Hodge structure on $H^1(A;\Q)$ be of
CM-type; this is the conventional definition for
$A$ to be of CM-type. This already implies that
the rational Hodge structure on $H^k(A;\Q)$ is of CM-type
for all $k=0,1,\cdots, 2n$. So, the abelian variety $A$
is of strong CM-type.
Let $({\cal K}, \Phi)$ be a CM pair of the abelian variety $A$;
${\cal K}$ is a commutative subalgebra of the Hodge endomorphism 
algebra of $H^1(A;\Q)$ such that $\dim_\Q {\cal K} = 2n$.
Then the reflex field $K^r$ of the CM pair $({\cal K},\Phi)$
is the Hodge endomorphism algebra of the level-$n$ subspace 
$[H^n(A;\Q)]_{\ell = n}$ (cf.\ \cite[Prop.\ 1.9.2]{MR257031}). 
% Shimura 70 bounded domain reflex for level n. 
So far, this is a well-known story. 
%  
% % CM-type \cite{borcea1998calabi}.
%  Borcea CY3
% 

In this note, we just assume that $A$ is of weak CM-type instead. 
The rational Hodge structure on the level-$n$ subspace $[H^n(A;\Q)]_{\ell = n}$
is assumed to be of CM-type; the Hodge endomorphism algebra 
$K':= {\rm End}([H^n(A;\Q)]_{\ell = n})^{\rm Hdg}$ is then a CM field; we call 
its degree $2n' := [K':\Q]$ the {\it reflex degree of} $A$. 
The following proof of 
Thms.\ \ref{thm:intro-main} and \ref{thm:intro-explct} begins
with the following observation: there exist a basis
$\{ dz^{a} \}_{a=1,\cdots,n}$ of $H^{1,0}(A;\C)$ and a basis
$\{ u^i\}_{i=1,\cdots, 2n}$ of $H^1(A;\Q)$ so that
\begin{align}
  \left( \begin{array}{c} dz^a \\  d\bar{z}^{\bar{a}} \end{array} \right)
  = \left( \begin{array}{cc} {\bf 1}_{n\times n} & \tau \\
    {\bf 1}_{n\times n} & \bar{\tau} \end{array} \right)
   \left( \begin{array}{c} u^1 \\ \vdots \\ u^{2n} \end{array} \right),  
   \qquad \qquad {\rm det}(\tau - \bar{\tau}) \neq 0,
     \label{eq:introduce-tau}
\end{align}
where $\tau$ is a complex-valued $n \times n$ matrix; 
the condition that
\begin{align*}
 \Omega = dz^1 \wedge \cdots \wedge dz^n
\end{align*}
is an eigenvector of the Hodge endomorphism field $K'$ implies that
the matrix $\tau$ takes value in the field $\phi'_{(n,0)}(K') \subset \C$,
where $\phi'_{(n,0)}$ is one of the $2n'$ embeddings of $K'$ into $\C$
given by assigning to $x \in K'$ its eigenvalue on $\Omega$.

%%%%%%%%%%%%%%%%%%%%%%%%%%%%%%%%%%%%%%%%%%%%%%
\subsection{The Reflex Degree 2}
%%%%%%%%%%%%%%%%%%%%%%%%%%%%%%%%%%%%%%%%%%%%%%

Shioda and Mitani showed \cite{shioda1974singular}
% Shioda Mitani 74
the following: for any weight-2 integral Hodge structure on
$H^2(T^4;\Z)$ with $h^{2,0}=1$ whose level-2 subspace is CM by
an imaginary quadratic field $K'$, there is a complex structure $I$
on $T^4$ so that the integral Hodge structure is reproduced
\cite{shioda1978period}, 
% Shioda Torelli for abelian surface
and the complex torus $(T^4,I)$ is isogenous to $E \times E$, where
$E$ is an elliptic curve with CM by $K'$. The statements
in Thms.\ \ref{thm:intro-main} and \ref{thm:intro-explct} with $n'=1$
are along the line of the result of Shioda--Mitani with a general $n$ not
necessarily $n=2$, but in a weaker form in that our question starts
off with an existing complex torus $(T^{2n},I)$ and its weight-$n$
Hodge structure than with an abstract weight-$n$ Hodge structure.
%; we also deal only with {\it rational} Hodge structures instead of {\it integral} Hodge structures.  

{\bf An elementary proof for the case $2n'=2$:}
let $\phi'_{(n,0)}(K')$ be $\Q(\sqrt{p})$,
where $p \in \Q_{<0}$, and $\sqrt{p}$ a purely imaginary complex number
with ${\rm Im}(\sqrt{p})>0$. The $\Q(\sqrt{p})$-valued matrix $\tau$
is expanded in the form of
\begin{align*}
 \tau = B_1 + B_2 \sqrt{p}, 
\end{align*}
where $B_1, B_2 \in M_n(\Q)$ ($\Q$-valued $n \times n$ matrix).
Then there must be $P, S \in M_n(\Q)$ such
that $P^{-1} B_2 S$ is diagonal and $\Q$-valued. Now,
\begin{align*}
  P^{-1} \left( \begin{array}{c} dz^1\\ \vdots \\ dz^n \end{array} \right)
   & \; = P^{-1}_{n\times n} \left( {\bf 1}_{n\times n}, \; (B_1+B_2\sqrt{p})_{n\times n}\right)
   \left( \begin{array}{cc} P & - B_1S \\ 0 & S \end{array} \right)
   \left( \begin{array}{cc} P^{-1} & P^{-1}B_1 \\ 0 & S^{-1} \end{array}
   \right) \left( \begin{array}{c} u^1 \\ \vdots \\ u^{2n} \end{array} \right)
   \nonumber \\
   & \; = \left( {\bf 1}_{n\times n}, \; ({\rm diag}) \sqrt{p} \right)
   \left( \begin{array}{cc} P^{-1} & P^{-1}B_1 \\ 0 & S^{-1} \end{array}
   \right) \left( \begin{array}{c} u^1 \\ \vdots \\ u^{2n} \end{array} \right).
\end{align*}
This implies that the complex torus in question is isogenous to
the product $\prod_{a=1}^n E_a$ of $n$ elliptic curves $\{E_a \}$
whose holomorphic 1-forms are $dw^{a}$ with $a=1,\cdots, n$ in 
\begin{align*}
  \left( \begin{array}{c} dw^{1} \\ \vdots \\ dw^{n} \end{array} \right) :=
  P^{-1} \left( \begin{array}{c} dz^1\\ \vdots \\ dz^n \end{array} \right),
   \qquad 
    \left( \begin{array}{c} u^{'1} \\ \vdots \\ u^{'2n}\end{array} \right) :=  
   \left( \begin{array}{cc} P^{-1} & P^{-1}B_1 \\ 0 & S^{-1} \end{array}
   \right) \left( \begin{array}{c} u^1 \\ \vdots \\ u^{2n} \end{array} \right), 
\end{align*}
and whose $H^1(E_a;\Q)$ are ${\rm Span}_\Q\{ u^{'a}, u^{'n+a} \}$.
Each of $E_a$'s is CM by $\Q(\sqrt{p})$. \qed

%%%%%%%%%%%%%%%%%%%%%%%%%%%%%%%%%%%%%%%%%%%%%%%
\subsection{The Reflex Degree 4}
%%%%%%%%%%%%%%%%%%%%%%%%%%%%%%%%%%%%%%%%%%%%%%%

Let us now prove Thms.\ \ref{thm:intro-main} and \ref{thm:intro-explct}
in the case the reflex degree is $2n'=4$. A CM field $K'$ of degree-4 
is classified by the Galois group ${\rm Gal}((K')^{\rm nc}/\Q)$, where $(K')^\mathrm{nc}$ is the normal closure of $\phi'_{(n,0)}(K')$ in $\C$.
There are three cases \cite[pp.\ 64--65, Ex.\ 8.4.(2)]{shimura2016abelian}:
% Shimura book CM abelian 
%
%
\begin{itemize}
\item [(A)] $K'/\Q$ is Galois, and 
${\rm Gal}(K'/\Q) \cong \Z_2\times \Z_2$. The field $K'$ is  
isomorphic to $\Q[x_1, x_2]/(x_1^2-p_1, \; x_2^2-p_2)$
for some negative rational numbers $p_1,p_2$ such that 
$\Q(\sqrt{p_1}) \neq \Q(\sqrt{p_2})$. 
\item [(B)] $K'/\Q$ is Galois, and 
${\rm Gal}(K'/\Q) \cong \Z_4$. The field $K'$ is isomorphic to 
$\Q[x,y]/(y^2-d, x^2-p-q y)$ for some positive square-free 
integer $d >1$ and rational numbers $p<0$ and $q$ such that 
$d' := p^2-q^2d \in d (\Q^\times)^2$. 
\item [(C)] $K'/\Q$ is not Galois, and 
${\rm Gal}((K')^{\rm nc}/\Q) \cong \Z_4 \rtimes \Z_2$. 
The field $K'$ is isomorphic to $\Q[x,y]/(y^2-d, \; x^2-p-q\eta)$ 
for some positive square free integer $d >1$ and rational numbers $p<0$ 
and $q$ such that $d' := p^2-q^2d \nin d (\Q^\times)^2$. 
\end{itemize}
The proof is given separately for the case (A), and for the case (B, C);
notations above will be used in the following. 

%%%%%%%%%%%%%%%%%%%%%%%%%%%%%%%%%%%%%%%
\subsubsection{The Case (A)}
\label{sssec:Rdeg4-caseA}
%%%%%%%%%%%%%%%%%%%%%%%%%%%%%%%%%%%%%%

{\bf proof:} Let us first prepare a little more notations. 
Let the generators $x_1$ and $x_2$ of the Hodge endomorphism algebra 
$K'$ be such that they are mapped by $\phi'_{(n,0)}$ to $\sqrt{p_1}$ and 
$\sqrt{p_2}$ (pure imaginary complex numbers in the upper half plane), 
respectively. The Galois group $\Z_2\vev{\sigma_1} \times \Z_2\vev{\sigma_2}$
is generated by $\sigma_1$ and $\sigma_2$, which acts on the Galois 
closure $(K')^{\rm nc} = \Q(\sqrt{p_1},\sqrt{p_2})$, as 
\begin{align*}
 \sigma_1: \sqrt{p_1} \longmapsto - \sqrt{p_1}, \qquad
           \sqrt{p_2} \longmapsto \sqrt{p_2}, \\
 \sigma_2: \sqrt{p_1} \longmapsto \sqrt{p_1}, \qquad 
          \sqrt{p_2} \longmapsto - \sqrt{p_2}. 
\end{align*}

The complex structure of an abelian variety $A$ under consideration 
is specified completely by the matrix $\tau$. It is only a necessary condition, though,
that the matrix $\tau$ is $\phi'_{(n,0)}(K')$-valued for the level-$n$ 
subspace to be of CM-type by $K'$. To be a sufficient 
condition, this $\phi'_{(n,0)}(K')$-valuedness condition has to be 
accompanied by one more: the Galois conjugate $\sigma_1(\Omega)$ 
in $H^n(A;\Q) \otimes \overline{\Q}$ has to be of pure Hodge $(p,n-p)$ type 
for some $0 < p < n$. 
%
% ; $\sigma_2(\Omega)$ and $\sigma_1 \circ \sigma_2(\Omega)$
% are of Hodge $(n-p,p)$ type and Hodge $(0,n)$ type, automatically. 
%

For $k=0,1, \cdots,n$, we call the following for 
$\Sigma \in H^n(A;\Q) \otimes \overline{\Q}$ the (\#$k$) condition:
\begin{itemize}
 \item [ (\#$k$) ] $\Sigma$ is orthogonal against any $n$-forms 
of the form $dz^{i_1} \wedge \cdots \wedge dz^{i_{n-k}}\wedge 
 d\bar{z}^{j_1} \wedge \cdots \wedge d\bar{z}^{j_k}$, where 
$\{ i_1,\cdots, i_{n-k} \}$ and $\{j_1,\cdots, j_k\}$ are subsets 
of $\{1,\cdots, n\}$.  
\end{itemize}
Therefore, the additional condition on the $\phi'_{(n,0)}(K')$-valued matrix $\tau$ 
is that $\Sigma = \sigma_1(\Omega)$ satisfies all the 
(\#$k$) conditions for $k \in \{ 0,\cdots, n\}$ except the 
(\#$p$) condition, and that the (\#$p$) condition is not satisfied 
for some choice $\{i_1,\cdots, i_{n-p}\}$ and $\{j_1,\cdots, j_p\}$. 

Here, we prepare a little more notations. Let $\tau^i$ denote the 
$i$-th row of the matrix $\tau$ (so, $\tau^i$ for each $i$ is in the 
$n$-dimensional vector space over $(K')^{\rm nc}$); 
here are other $n$-component vectors:
\begin{align*}
 \delta^i := \tau^i - (\tau^i)^{\sigma_1}, \qquad 
 \epsilon^i := (\tau^i)^{\sigma_1\sigma_2} - (\tau^i)^{\sigma_1} =
    - \left(\tau^i - (\tau^i)^{\sigma_2}\right)^{\sigma_1}, \qquad 
 (\epsilon^i)^{\sigma_1\sigma_2} = \tau^i - (\tau^i)^{\sigma_2}. 
\end{align*}
The vector $\sqrt{p_1}\delta^i$ is in the $n$-dimensional vector space over
$\Q(\sqrt{p_2})$, and $\sqrt{p_2} \epsilon^i$ is in the $n$-dimensional
vector space over $\Q(\sqrt{p_1})$, in fact.  

With these notations, the condition (\#$k$) on $\Sigma = \sigma_1(\Omega)$ 
is equivalent to 
\begin{align*}
 {\rm det}\left( \delta^{i_1T}, \cdots, \delta^{i_{n-k}T}, \epsilon^{j_1T},\cdots 
\epsilon^{j_kT} \right) = 0; 
\end{align*}
that is, those $\Q(\sqrt{p_1},\sqrt{p_2})$-valued $n$ vectors are not linearly 
independent. The additional condition referred to above is equivalent to 
the condition that 
\begin{align}
{\rm Span}_{\C}\{ \delta^1, \cdots, \delta^n\} \cong \C^{n-p}, \quad 
{\rm Span}_{\C} \{ \epsilon^1, \cdots, \epsilon^n \} \cong \C^p, \quad
{\rm Span}_{\C} \{ \delta's, \epsilon's \} \cong \C^n;
\label{eq:deltaNeps-spanAllTgthr}
\end{align}
the field $\C$ above may have been replaced by $\Q(\sqrt{p_1},\sqrt{p_2})$. 

One may extract a little more information on which combination 
of $\delta$'s and $\epsilon$'s are linearly independent and generate 
the whole $n$-dimensional vector space. Noting that 
${\rm det}(\tau - \bar{\tau}) \neq 0$ and 
$\tau - \bar{\tau} = \delta - \epsilon$, where 
$\delta := \tau - (\tau)^{\sigma_1}$ and 
$\epsilon := \tau^{\sigma_1\sigma_2}- \tau^{\sigma_1}$ are $n \times n$ matrices, 
one can conclude that $\{ \delta^{i_1},\cdots, \delta^{i_{n-p}}, \epsilon^{j_1},\cdots, \epsilon^{j_p}\}$ is linearly independent for at least one 
mutually exclusive choice of $\{i_1,\cdots, i_{n-p}\} \amalg 
\{j_1,\cdots, j_p\} = \{ 1,\cdots, n\}$. Let us rename 
$\{u^1,\cdots, u^n\}$ and $\{dz^1, \cdots, dz^n\}$ so that 
$\{ \delta^1,\cdots, \delta^{n-p}, \epsilon^{n-p+1},\cdots, \epsilon^{n} \}$ 
are linearly independent.  

There must be $\Q(\sqrt{p_1})$-valued $(n-p) \times p$ matrix $c_1$
and $\Q(\sqrt{p_2})$-valued $p \times (n-p)$ matrix $c_2$ such that 
\begin{align*}
  \epsilon^a + (c_1)^{ab} \epsilon^{n-p+b} & \; = 0, \qquad
      a=1,\cdots, n-p, \quad b = 1,\cdots, p, \\
  (c_2)^{bd} \delta^d + \delta^{n-p+b} & \; = 0, \qquad
       b = 1,\cdots, p, \quad d = 1,\cdots, n-p,  
\end{align*}
because all of $\sqrt{p_1}^{-1} \delta^i$'s are $\Q(\sqrt{p_2})$-valued 
$n$-component vectors, and all of $\sqrt{p_2}^{-1} \epsilon^i$'s are 
$\Q(\sqrt{p_1})$-valued $n$-component vectors. We may use the two matrices 
to redefine $n$ holomorphic coordinates on $A$ by\footnote{
It is possible to prove that this redefinition is invertible; we use 
${\rm det}(\tau - \bar{\tau}) \neq 0$ for this.  
} % 
\begin{align*}
 dw^a & \; = dz^a + ((c_1)^{\sigma_1})^{ab} dz^{n-p+b}, \\
 dw^{n-p+b} & \; = (c_2)^{bd} dz^d+ dz^{n-p+b}; 
\end{align*}
this means that
\begin{align*}
 \left( \begin{array}{c} dw^a|_{a=1,\cdots, n-p} \\ dw^{n-p+b}|_{b=1,\cdots, p}
   \end{array} \right) =
 \left( \begin{array}{cc|c}
    {\bf 1}_{(n-p)\times (n-p)} &  (c_1)^{\sigma_1} & (1, (c_1)^{\sigma_1}) \tau \\
    (c_2) & {\bf 1}_{p\times p} & ((c_2),1) \tau \end{array} \right) 
  \left( \begin{array}{c} u^1 \\ \vdots \\ u^{2n} \end{array} \right).
\end{align*}
The top $(n-p)$ rows of the $n \times 2n$ matrix above
are $\Q(\sqrt{p_1})$-valued and the bottom $p$ rows $\Q(\sqrt{p_2})$-valued,
as one may notice that 
\begin{align*}
  \tau  & \; = B_1 + B_3 \sqrt{p_2} + 2^{-1} \delta
     = B_1 + B_2 \sqrt{p_1} - 2^{-1} \epsilon^{\sigma_1},  
\end{align*}
when we expand the matrix $\tau$ as 
\begin{align*}
  \tau & \; = B_1 + B_2 \sqrt{p_1} + B_3 \sqrt{p_2} + B_4 \sqrt{p_1}\sqrt{p_2},
  \qquad B_1, B_2, B_3, B_4 \in M_n(\Q). 
\end{align*}
It is not hard to find a pair of invertible matrices
$P^{-1} \in M_{n-p}(\Q(\sqrt{p_1})) \oplus M_{p}(\Q(\sqrt{p_2}))$ and
$S\in M_{2n}(\Q)$ so that the $n \times 2n$ matrix above is multiplied
by $P^{-1}$ from the left and by $S$ from the right to be of the form
\begin{align*}
  \left( \begin{array}{cc|cc}
    {\bf 1}_{(n-p)\times (n-p)} & & \sqrt{p_1} {\bf 1}_{(n-p)\times (n-p)} & \\
    & {\bf 1}_{p\times p} & & \sqrt{p_2} {\bf 1}_{p\times p} \end{array} \right). 
\end{align*}
This indicates that the abelian variety $A$ under consideration
is isogenous to $(E_{\sqrt{p_1}})^{n-p} \times (E_{\sqrt{p_2}})^p$, where
$E_{\sqrt{p_1}}$ is an elliptic curve with CM by $\Q(\sqrt{p_1})$ and
$E_{\sqrt{p_2}}$ an elliptic curve with CM by $\Q(\sqrt{p_2})$. \qed

\begin{props}
 {\it  Let $A$ be a weak CM-type abelian variety of $n$-dimensions, and
  the CM field $K'$ on the level-$n$ subspace is of degree-$(2n'=4)$ and
  is in the case (A) in the classification
  of} \cite[pp.\ 64--65, Ex.\ 8.4.(2)]{shimura2016abelian}.
  {\it Then the level-$n$ subspace consists of Hodge components
  with $h^{n,0} = h^{0,n}=1$ and $h^{p,n-p} = h^{n-p,n}=1$ when
  $A$ is isogenous to $E_1^{n-p} \times E_2^{p}$ for $0 < p < n$ and
  some mutually non-isogenous CM elliptic curves $E_1$ and $E_2$. }
   \qed
\end{props}

%%%%%%%%%%%%%%%%%%%%%%%%%%%%%%%%%%%%%%%%
\subsubsection{The Case (B, C)}
\label{sssec:Rdeg4-caseBC}
%%%%%%%%%%%%%%%%%%%%%%%%%%%%%%%%%%%%%%%%

{\bf proof:} Let us prepare a little more notations.
For a CM field $K'$ in the case (B), let the generators
$y$ and $x$ of $K'$ be such that they are mapped
by $\phi'_{(n,0)}$ to a real positive number $\sqrt{d}$
and a pure imaginary number in the upper half plane
$\xi_+ := \sqrt{p+q\sqrt{d}}$, respectively. The Galois closure
$(K')^{\rm nc}$ of $\phi'_{(n,0)}(K') \subset \C$ is
$\Q (\sqrt{d},\xi_+) = \phi'_{(n,0)}(K')$ itself. 
An order-4 generator $\sigma_0$ of the Galois group
$\Z_4 \cong {\rm Gal}((K')^{\rm nc}/\Q)$ is
\begin{align*}
  \sigma_0: \sqrt{d} \longmapsto - \sqrt{d}, \quad
  \xi_+ = \sqrt{p+q\sqrt{d}} \longmapsto \sqrt{p-q\sqrt{d}} =: \xi_-
    = - \sqrt{d'}/\xi_+,
\end{align*}
%
% which also means that
% %
% \begin{align}
%  \sqrt{p-q\sqrt{d}} = - \frac{\sqrt{d'}}{\sqrt{p+q\sqrt{d}}} \longmapsto
%   \frac{\sqrt{d'}}{\sqrt{p-q\sqrt{d}}} = - \sqrt{p+q\sqrt{d}} . 
% \end{align}
%
where $\sqrt{d'}>0$.
The complex conjugation is $\sigma_0^2 \in \Z_4$. 

For a CM field $K'$ in the case (C), let the generators
$\eta$ and $\xi$ of $K'$ be such that they are mapped by $\phi'_{(n,0)}$
to the real positive $\sqrt{d}$ and the pure imaginary $\xi_+$ in the upper
half plane, respectively. The Galois
closure $(K')^{\rm nc}$ of $\phi'_{(n,0)}(K') \subset \C$ is
\begin{align*}
 &  {\rm Span}_\Q\left\{ 1, \sqrt{d}, \sqrt{d'}, \sqrt{d}\sqrt{d'},
     \xi_{\pm}, \sqrt{d}\xi_{\pm} 
     \right\} \\
     & \; \cong \Q[y,z,x_+,x_-]/
          \left( y^2-d, z^2-d', x_\pm^2 - p \pm qy, x_+x_-+z \right). 
\end{align*}
The Galois group ${\rm Gal}((K')^{\rm nc}/\Q)$ is isomorphic to\footnote{
i.e., $\sigma_3 \circ \sigma_0 \circ \sigma_3 =  \sigma_0^3 = \sigma_0^{-1}$
} %
$\Z_4\vev{\sigma_0} \rtimes \Z_2\vev{\sigma_3}$, where $\sigma_0$ is
the Galois transformation 
\begin{align*}
  \sqrt{d} \longmapsto -\sqrt{d}, \quad \sqrt{d'} \longmapsto - \sqrt{d'},
  \quad \!\! 
  \xi_+ \longmapsto \xi_-, \quad \xi_- \longmapsto - \xi_+, 
\end{align*}
and $\sigma_3$ the transformation 
\begin{align*}
  \sqrt{d} \longmapsto - \sqrt{d}, \quad \sqrt{d'} \longmapsto \sqrt{d'},
  \quad
  \xi_+ \longleftrightarrow \xi_-. 
\end{align*}
The complex conjugation is $\sigma_0^2$; the Galois transformation
$\sigma_0^3 \circ \sigma_3$ on $K^{'{\rm nc}}$ is trivial when restricted
to the subfield $\phi'_{(n,0)}(K') = \Q(\sqrt{d}, \xi_+)$.
The four embeddings of $K'$ are $\phi'_{(n,0)}$, $\sigma_0 \circ \phi'_{(n,0)}$,
$\sigma_0^2 \circ \phi'_{(n,0)}$ and $\sigma_0^3 \circ \phi'_{(n,0)}$. 

Now, a proof starts here. Just as in the case studied in
section \ref{sssec:Rdeg4-caseA}, it is necessary for the $\tau$
matrix in (\ref{eq:introduce-tau}) to be $\phi'_{(n,0)}(K')$-valued
so that the level-$n$ subspace of an abelian variety under consideration
is CM by $K'$. For a necessary and sufficient condition, just one
more condition has to be added: the Galois conjugate $\Sigma :=
\sigma_0(\Omega)$ is of pure Hodge $(p,n-p)$ type for some $0 < p < n$. 
This condition is also translated into the language of the (\#$k$)
conditions, as we have done in section \ref{sssec:Rdeg4-caseA}.

Here, we prepare a little more notations. Each row vector $\tau^i$
(the $i$-th row) of the matrix $\tau$ is in the $n$-dimensional
vector space over $(K')^\mathrm{nc}$. Here are other $n$-component
vectors: 
\begin{align*}
  \delta^i := \tau^i - (\tau^i)^{\sigma_0}, \qquad
  \epsilon^i = (\tau^i)^{\sigma_0^2} - (\tau^i)^{\sigma_0}. 
\end{align*}
When the CM field $K'$ is in the case (B, C), there is a relation
$\epsilon^i = - \sigma_0(\delta^i)$. 

The relation $\epsilon^i = - \sigma_0(\delta^i)$ implies that
the subspace spanned by $\{ \delta^i\}$ has the same dimension
as that spanned by $\{ \epsilon^i \}$. We should be in a situation
(\ref{eq:deltaNeps-spanAllTgthr}), but with $p = n-p$.
It follows immediately that there is no such things as a
weak CM abelian variety $A$ with a reflex degree-4 
CM field $K'$ in the case (B, C), if $\dim_\C A = n$ is odd. 
In the rest of this section \ref{sssec:Rdeg4-caseBC},
$n=: 2r$ for an integer $r > 0$. 

We may think that the $r$ vectors in $\C^n$,
$\{ \delta^{r+1}, \cdots, \delta^{n}\}$, are linearly independent over $\C$,
by allowing ourselves to rename $dz^1, \cdots, dz^n$.
There then exists an $r \times r$ matrix $M$ such that
\begin{align*}
  \left( \begin{array}{c} \delta^1 \\ \vdots \\ \delta^r \end{array} \right)
  = M  \left( \begin{array}{c} \delta^{r+1} \\ \vdots \\ \delta^n \end{array}
  \right); 
\end{align*}
it then follows that
\begin{align*}
   \left( \begin{array}{c} \epsilon^1 \\ \vdots \\ \epsilon^r \end{array} \right)
  = M^{\sigma_0}  \left( \begin{array}{c} \epsilon^{r+1} \\ \vdots \\ \epsilon^n \end{array}
  \right). 
\end{align*}
The matrix $M$ is not just $\C$-valued, but takes value in fact in
the subfield
\begin{align*}
  (K')^r := {\rm Span}_\Q \left\{ 1, \; \sqrt{d'}, \;
   \left( \xi_+ + \xi_-\right), \; \sqrt{d}(\xi_+-\xi_-)
      \right\} 
\end{align*}
of $(K^{'{\rm nc}}) \subset \C$. To see this, it is enough to use
an expansion
\begin{align*}
  \tau = B_1 + B_2 \sqrt{d} + B_3 \xi_+ + B_4 \sqrt{d} \xi_+, \qquad
     B_1, B_2, B_3, B_4 \in M_n(\Q), 
\end{align*}
and do the computation:
\begin{align*}
  \delta & \; = \sqrt{d} \left[ 2B_2 + B_3 (\xi_+ - \xi_-)/\sqrt{d}
    + B_4 (\xi_+ + \xi_-) \right], \\
  \epsilon & \; = \sqrt{d} \left[ 2B_2 +B_3(\xi_+ - \xi_-)/\sqrt{d}
    + B_4 (\xi_+ + \xi_-) \right]^{\sigma_0}. 
\end{align*}
The subfield $(K')^r$ above is the reflex field of the pair
of the CM field $K'$ and its half a set of embeddings
$\{ \phi'_{(n,0)}, \sigma_0 \circ \phi'_{(n,0)}\}$.

So far, the relation (\ref{eq:introduce-tau}) remains in its form,
when we allow ourselves to rename $u^1, \cdots, u^n$ if necessary
(as we have renamed $\{ dz^1, \cdots, dz^n\}$). Now, 
we may redefine the $n$ holomorphic coordinates on $A$ as follows:
\begin{align*}
 \left( \begin{array}{c} dw^a \\ dw^{r+b} \end{array} \right)
  := \left( \begin{array}{cc} \mathbf{1}_{r\times r} & -M \\ \mathbf{1}_{r\times r} & -M^{\sigma_0^3}
    \end{array} \right)
    \left( \begin{array}{c} dz^1 \\ \vdots \\ dz^n \end{array} \right). 
\end{align*}
Using the fact that
\begin{align*}
  \tau - B_1 &\; = \frac{1}{2} \delta
  - \frac{1}{4}\epsilon + \frac{1}{4} \epsilon^{\sigma_0^2}
  = \left( \begin{array}{c} 
    \frac{1}{2} M \underline{\delta}
    - \frac{1}{4} M^{\sigma_0} \underline{\epsilon}
    + \frac{1}{4} M^{\sigma_0^3} (\underline{\epsilon})^{\sigma_0^2} \\
    \frac{1}{2} \underline{\delta}
    - \frac{1}{4} \underline{\epsilon}
    + \frac{1}{4} (\underline{\epsilon})^{\sigma_0^2}
  \end{array} \right) \\
  & \; =
  \frac{1}{4} \delta - \frac{1}{4} \delta^{\sigma_0^2}
  + \frac{1}{2}\epsilon^{\sigma_0^2} =
  \left( \begin{array}{c}
    \frac{1}{4} M \underline{\delta}
    - \frac{1}{4} M^{\sigma_0^2} (\underline{\delta})^{\sigma_0^2}
    + \frac{1}{2} M^{\sigma_0^3} (\underline{\epsilon})^{\sigma_0^2} \\
    \frac{1}{4} \underline{\delta}
    - \frac{1}{4} (\underline{\delta})^{\sigma_0^2}
    + \frac{1}{2} (\underline{\epsilon})^{\sigma_0^2}
     \end{array} \right), 
\end{align*}
where $\underline{\delta}$ and $\underline{\epsilon}$ are the bottom
half $r\times n$ blocks of the $n \times n$ matrices $\delta$ and $\epsilon$,
respectively, we have
\begin{align}
 \left( \begin{array}{c} dw^a \\ dw^{r+b} \end{array} \right) = 
 \left( \begin{array}{cc|c} {\bf 1}_{r\times r} & - M &
   ({\bf 1},-M)B_1
   + \frac{M-M^{\sigma_0}}{4}\underline{\epsilon}
   - \frac{M-M^{\sigma_0^3}}{4}(\underline{\epsilon})^{\sigma_0^2} \\
   {\bf 1}_{r\times r} & - M^{\sigma_0^3} &
   ({\bf 1},-M^{\sigma_0^3})B_1 
   + \frac{M-M^{\sigma_0^3}}{4}\underline{\delta}
  - \frac{M^{\sigma_0^2}-M^{\sigma_0^3}}{4} (\underline{\delta})^{\sigma_0^2}
   \end{array} \right) \left( \begin{array}{c} u^1 \\ \vdots \\ u^{2n}
 \end{array} \right).
 \label{eq:prd-std-form-BC-pre}
\end{align}
The bottom $r$ rows of the $n \times 2n$ matrix above are 
the Galois conjugate by $\sigma_0^3$ of the top $r$ rows
(remember that $\epsilon^{\sigma_0^3} = -\delta$). 

The top $r$ rows of the $n \times 2n$ matrix above are 
$(K')^r$-valued, not just in the upper left $r \times n$ block
$({\bf 1}_{r \times r}, -M)$, but also in the upper right $r \times n$ block.
To see this, note first that the last two terms of the upper right
$r \times n$ block is  
\begin{align*}
 \frac{1}{4} (M-M^{\sigma_0})\underline{\epsilon}
  - \frac{1}{4}(M-M^{\sigma_0^3})(\underline{\epsilon})^{\sigma_0^2} = 
  - \frac{1}{4} \left( M^{\sigma_0}\underline{\epsilon}
   - (M^{\sigma_0}\underline{\epsilon})^{\sigma_0^2}\right)
   + \frac{1}{4} M \left( \underline{\epsilon}
   - (\underline{\epsilon})^{\sigma_0^2}\right).
\end{align*}
Both $M^{\sigma_0} \underline{\epsilon}$ and $\underline{\epsilon}$
in the 1st and 2nd term, respectively, on the right-hand side
take value in $\sqrt{d} \times \sigma_0((K')^r)$. Any element
$\xi \in \sigma_0((K')^r)$ has the property that
$\xi - \xi^{\sigma_0^2} \in \sqrt{d} (K')^r$. It is now easy to see
that the $r \times n$ matrix on the right hand side takes value
in $(K')^r$. 

It is not hard to find a pair of invertible matrices
$P^{-1} \in M_{r}((K')^r)$ and $S \in M_{2n}(\Q)$ so that the
upper $r \times 2n$ block of the $(K')^r$-valued matrix above
is multiplied by $P^{-1}$ from the left and by $S$ from the right
to be in the form of
\begin{align}
  \left( \begin{array}{cccc}
   {\bf 1}_{r \times r}, & \; \sqrt{d'} {\bf 1}_{r \times r},
    & \; (\xi_++\xi_-) {\bf 1}_{r \times 1}, & \;
          \sqrt{d'}(\xi_+-\xi_-){\bf 1}_{r \times r}
  \end{array} \right).
  \label{eq:prd-std-form-BC}
\end{align}
We may multiply the matrix $(P^{-1})^{\sigma_0^3}$ to redefine the coordinates
$dw^{r+b}$'s, to convert the bottom $r \times 2n$ block of the
$\sigma_0^3((K')^r)$-valued matrix in (\ref{eq:prd-std-form-BC-pre})
to the $\sigma_0^3$-transform of (\ref{eq:prd-std-form-BC}). 
This indicates that the abelian variety $A$ under consideration
is isogenous to the product of $r$ copies of an abelian surface
with CM by $(K')^r$. \qed 

\begin{props}
\label{props:caseBC-level2-HdgStr}
{\it Let $A$ be a weak CM-type abelian variety of $n$-dimensions. Suppose further that 
the CM field $K'$ on the level-$n$ subspace is of degree-$(2n'=4)$, and 
is in the case (B, C) in the classification of }
\cite[pp.\ 64--65, Ex.\ 8.4.(2)]{shimura2016abelian}. 
{\it Then the level-$n$ subspace consists of Hodge components with 
$h^{n,0}=h^{0,n}=1$ and $h^{r,r}=2$, where $n=2r$. }    \qed
\end{props}

%%%%%%%%%%%%%%%%%%%%%%%%%%%%%%%%%%%%%%%%%%%%%%%
\section{K3$\times T^2$ of Weak CM-type}
\label{sec:K3xT2}
%%%%%%%%%%%%%%%%%%%%%%%%%%%%%%%%%%%%%%%%%%%%%%%

Let $S$ and $E$ be a complex K3 surface and an elliptic curve over $\C$, 
respectively, and $X=S \times E$ in this section \ref{sec:K3xT2}. 
We assume that $X$ is of weak CM-type, which means that the 
level-$(n=3)$ subspace of $H^3(X;\Q) \cong H^2(S;\Q) \otimes_\Q H^1(E;\Q)$
has a CM-type rational Hodge structure, and read out implications 
of this assumption. Let $K'$ be the Hodge endomorphism algebra 
of the level-3 subspace; $2n' := [K':\Q]$.  

{\bf proof of Thm.\ \ref{thm:K3xT2}:} Set a holomorphic (2,0)-form on $S$ and a holomorphic (1,0)-form on $E$ as follows:
\begin{align*}
 \Omega_S = \sum_{i=1}^{22} e_i \Pi_i, \qquad 
 \Omega_E = \hat{\alpha} + \tau_E \hat{\beta},
\end{align*}
where $\{e_{i=1,\cdots, 22} \}$ and $\{ \hat{\alpha}, \hat{\beta}\}$ 
are bases of $H^2(S;\Q)$ and $H^1(E;\Q)$, respectively.  We choose 
the normalization of $\Omega_S$ such that at least one of $\Pi_i$'s is in $\Q$. 
It is then a necessary condition for $S \times E$ to be of weak CM-type that 
\begin{align*}
  \Q(\Pi_i,\tau_E)=\phi'_{(3,0)}(K'); 
\end{align*}
$\phi'_{(3,0)}$ is the embedding $K' \hookrightarrow \C$ given by the 
eigenvalue of the endomorphisms on the Hodge (3,0)-type $\Omega_S \Omega_E$. 
To make it a sufficient condition as well, the following condition has 
to be added: when $G := {\rm Gal}((K')^{\rm nc}/\Q)$ and $H$ is the subgroup 
corresponding to the subfield $\phi'_{(3,0)}(K')$ within its normal 
closure $(K')^{\rm nc} \subset \C$, there must be $(n'-1)$ $H$-cosets 
of the form $\sigma H$ (among the $2n'$ cosets in $G/H$) for which 
$\Omega_S^\sigma\Omega_E^\sigma$ is a (2,1)-form. 
%The set of such cosets (embeddings of $K'$ into $(K')^\mathrm{nc}$) will be denoted by $\Phi^{(2,1)}$. 

%The conditions that $\sigma \circ \phi'_{(n,0)}(\Omega_E\Omega_S)$ is 
%free from the Hodge-(3,0) component, and from the Hodge-(0,3) %component, 
%are given by 
%
%\begin{align}
% \langle \Omega_E^\sigma, \; \overline{\Omega}_E \rangle
% (\Omega_S^\sigma, \; \overline{\Omega}_S )  = 0,  \qquad 
% \langle \Omega_E^\sigma, \; \Omega_E \rangle 
% (\Omega_S^\sigma, \; \Omega_S) =0; 
%  \label{eq:K3xT2-cond-A}
%\end{align}
%
%for $\sigma \circ \phi'_{(n,0)}(\Omega_E\Omega_S)$ to be free from 
%Hodge-(1,2) component, we should have 
%
%\begin{align}
% \langle \Omega_E^\sigma, \; \overline{\Omega}_E \rangle 
% ( \Omega_S^\sigma, \; \Omega_S )  = 0, \qquad 
% \langle \Omega_E^\sigma, \; \Omega_E \rangle 
% ( \Omega_S^\sigma, \omega_{a} )= 0, \quad a=1,\cdots, 20, 
%   \label{eq:K3xT2-cond-B}
%\end{align}
%
%where $\omega_a$'s are chosen to span $H^{1,1}(S;\R)$ of 20 dimensions. 
%Noting that it is impossible for $\Omega_E^\sigma$ to be both orthogonal 
%to $\Omega_E$ and $\overline{\Omega}_E$, one finds that such a $\sigma$ 
%satisfying both (\ref{eq:K3xT2-cond-A}, \ref{eq:K3xT2-cond-B}) 
%is either one of the two types: 
%

Note that the 2-form $\Omega_S^\sigma$ on $S$ and the 1-form $\Omega_E^\sigma$ on $E$ are either one of the following two types for their product to be a (2,1)-form.
\begin{align*}
 (*1): & \qquad \qquad \quad
  \Omega_S^\sigma = \Omega_S, \qquad \quad 
%  \langle \Omega_E^\sigma , \overline{\Omega}_E \rangle = 0, \qquad 
  \Omega_E^\sigma = \overline{\Omega}_E,  \\
 (*2): & \qquad \qquad \quad
   \Omega_S^\sigma \in H^{1,1}(S;\C), \qquad \quad 
%   \langle \Omega_E^\sigma, \Omega_E \rangle = 0, \qquad 
  \Omega_E^\sigma = \Omega_E. 
\end{align*}
Now, the Galois orbit of $\tau_E \in (K')^{\rm nc}$  consists of just 
$\{ \tau_E, \overline{\tau}_E \}$;
when $\sigma^{(1)}$ is in type (*1) and $\sigma^{(2)}$ in type (*2),
the Galois images of $\tau_E$ under 1, $\rho\sigma^{(1)}$ and $\sigma^{(2)}$ are $\tau_E$, while the images under $\rho$, $\sigma^{(1)}$ and $\rho\sigma^{(2)}$ are $\overline{\tau}_E$, where $\rho$ is the complex conjugate.
So the minimal polynomial of $\tau_E$ over $\Q$ is degree-2, and 
$\Q(\tau_E)$ is an imaginary quadratic field. $E$ is CM, in particular.

The transcendental lattice of the K3 surface $S$ is also of CM-type. 
This is because all the Galois conjugates of $\Omega_S$ are  
$\Omega_S$ itself (under 1 and $\sigma^{(1)}$), 
$\overline{\Omega}_S$ (under $\rho$ and $\rho \sigma^{(1)}$), and $\Omega_S^\sigma \in H^{1,1}(S;\C)$ (under $\sigma^{(2)}$ and $\rho\sigma^{(2)}$); all the Galois conjugates 
of the $K'$-eigenvector $\Omega_S$ are of some pure Hodge types.  
In particular, a weak CM $S \times E$ is also of strong CM-type. \qed

\begin{rmk}
There are two distinct situations covered in the proof above
(although the proof did not need separate treatments): 
one is when $\Q(\tau_E) \cap (\Q(\Pi_i))^{\rm nc} = \Q$ so that $[(K')^{\rm nc}:\Q] = [\Q(\tau_E):\Q][(\Q(\Pi_i))^{\rm nc}:\Q]$
and the other is when $\tau_E \in (\Q(\Pi_i))^{\rm nc}$ so that
$[(K')^{\rm nc}:\Q] = [(\Q(\Pi_i))^{\rm nc}:\Q]$. In the first situation, 
the level-3 subspace is $T_{S\Q} \otimes H^1(E;\Q)$, where $T_{S\Q}$ is 
the transcendental lattice of $S$ tensored by $\Q$. In the second situation, 
the level-3 subspace is a proper subspace of $T_{S\Q} \otimes H^1(E;\Q)$, 
half in the dimension over $\Q$ (cf.\ \cite{Kanno:2017nub}). 
\end{rmk}

\begin{rmk}
\label{rmk:dcmps}
The proof of Thm. \ref{thm:K3xT2} above also contains an implication 
for a variety $Y$ with a trivial canonical bundle much more general than 
$Y = {\rm K3} \times T^2$. When $X_1$ and $X_2$ are both projective 
varieties with trivial canonical bundles, and $Y = X_1 \times X_2$ is 
of weak CM-type with the Hodge endomorphism field $K'$, then one can 
repeat the argument in the proof to 
see for any $\sigma \in {\rm Gal}((K')^{\rm nc}/\Q)$ that 
both $\Omega_1^{\sigma}$ and $\Omega_2^{\sigma}$ are of pure Hodge type;
here $\Omega_1$ and $\Omega_2$ are non-zero sections of the 
trivial canonical bundles ${\rm det}(T^*X_1)$ and ${\rm det}(T^*X_2)$, 
respectively. This is enough for Thm. \ref{thm:dcmps}.

Thm. \ref{thm:K3xT2} can be regarded as a corollary of Thm. \ref{thm:dcmps}.
This is because both K3 and $T^2$ are of weak CM-type when 
$Y={\rm K3} \times T^2$ is of weak CM-type (Thm. \ref{thm:dcmps}), 
and weak CM K3 surfaces and elliptic curves are always of strong CM-type. 
\end{rmk}

%%%%%%%%%%%%%%%%%%%%%%%%%%%%%%%%%%%%%%%%%%%%%%
\section{Open Problems}
\label{sec:open-prblm}
%%%%%%%%%%%%%%%%%%%%%%%%%%%%%%%%%%%%%%%%%%%%%%

The present authors are not aware of a logic to infer that a weak 
CM-type complex projective non-singular variety $X$ with a trivial 
canonical bundle ${\rm det}(T^*X)$ is also of strong CM-type, except in 
the special cases covered by Thms.\ \ref{thm:intro-main}, 
\ref{thm:intro-explct} and \ref{thm:K3xT2}. On the other hand, not 
one explicit example of a weak CM-type $X$ that is not of strong CM-type 
is known to the authors. 
The expositions in this section \ref{sec:open-prblm} are intended to 
draw attention of mathematicians\footnote{
This section does not contain a new mathematical result. 
This material is still included in this arXiv 
preprint, as there might be a small number of readers who might  
find this useful.  
} % 
 to the question whether the weak CM
property also implies the strong CM property or not.

This question is trivial when $X$ is an elliptic curve or a K3 surface;
the Neron--Severi lattice within $H^2({\rm K3};\Z)$ has enough 
Hodge endomorphisms. This is a non-trivial question when $X$ is an abelian 
surface, but we have a complete answer 
(Thms.\ \ref{thm:intro-main} and \ref{thm:intro-explct}). So, the answer 
to this question is yes for $X$ with its dimension $n=1,2$. 

%%%%%%%%%%%%%%%%%%%%%%%%%%%%%%%%%%%%%%%%
\subsection{Abelian Varieties with a Higher Reflex Degree}
\label{ssec:open-abel}
%%%%%%%%%%%%%%%%%%%%%%%%%%%%%%%%%%%%%%%%

One may ask whether Thms.\ \ref{thm:intro-main} and \ref{thm:intro-explct}
still hold true for abelian varieties $A$ of arbitrary dimension $n$ 
and arbitrary reflex degree $2n'$, not just for $n'=1, 2$. If one 
is to approach this question through the case-by-case approach as 
in section \ref{sec:abel}, the language introduced by 
Dodson \cite{MR735406} (reviewed below) will be useful for the 
classification. 

\begin{notn}
For $N \in \N$, the group ${\rm Im}(N,2)$ is defined to be the subgroup 
of the permutation group $\mathfrak{S}_{2N}$ of $2N$ elements explained 
shortly. Introduce to the set of $2N$ elements a structure of $N$ mutually 
exclusive pairs ($\phi_i$ and $\bar{\phi}_i$ in $\{ \phi_1, \bar{\phi}_1, \phi_2, \bar{\phi}_2, \cdots, 
\phi_N, \bar{\phi}_N \}$ are regarded as a pair), and fix this additional information. 
The subgroup ${\rm Im}(N,2)$ consists of permutations 
$ \sigma \in \mathfrak{S}_{2N}$ that preserve the pairs; 
when $\sigma$ maps $\phi_i$ to either one of 
$\{ \phi_{[\sigma](i)}, \bar{\phi}_{[\sigma](i)}\}$, where 
$[\sigma] \in \mathfrak{S}_N$, then $\sigma(\bar{\phi}_{i})$ 
has to be the other one in $\{ \phi_{[\sigma](i)}, \bar{\phi}_{[\sigma](i)}\}$
for the permutation $\sigma \in \mathfrak{S}_{2N}$ to be in ${\rm Im}(N,2)$.   
There is a structure that 
\begin{align}
{\rm Im}(N,2) \cong (\Z_2)^N \rtimes \mathfrak{S}_N.
 \label{eq:ImN2-str}
\end{align}
\end{notn}

\begin{lemma}  \cite[\S 1--2]{MR735406}
For a given $N \in \N$, a subgroup $G$ of ${\rm Im}(N,2)$ is specified 
uniquely by a triple $(G_0, (\Z_2)^v, s)$ satisfying the following three conditions. 
\begin{itemize}
\item $G_0$ is a subgroup of $\mathfrak{S}_N$ that acts transitively 
on the $N$ elements.
\item $(\Z_2)^v$ is a subgroup of $(\Z_2)^N \subset {\rm Im}(N,2)$
for an integer $v$ in the range $0 < v \leq N$ where
\begin{itemize}
\item the $(\Z_2)^v$ 
subgroup of $(\Z_2)^N$ is preserved by the action of any $[g_0]\in G_0$;
$[g_0] \in G_0 \hookrightarrow \mathfrak{S}_N$ acts on $(\Z_2)^N \subset 
{\rm Im}(N,2)$ through the adjoint action (cf (\ref{eq:ImN2-str})), and 
\item the $(\Z_2)^v$ group contains the diagonal subgroup $\Z_2$ of 
$(\Z_2)^N$; the generator of the diagonal $\Z_2$ is denoted by $\rho$. 
\end{itemize}
\item $s \in Z^1(G_0; (\Z_2)^N/(\Z_2)^v)$.  
\end{itemize}
The subgroup $G$ of a triple, denoted by $G(G_0, (\Z_2)^v, s)$ or $G(G_0,v,s)$, 
is 
\begin{align}
 G = \bigcup_{g \in G_0} \bigcup_{x \in s(g)} (x,g) \subset (\Z_2)^N \rtimes \mathfrak{S}_N
  \cong {\rm Im}(N,2), 
\label{eq:G-explicit-formula}
\end{align}
where the $(\Z_2)^N/(\Z_2)^v$-valued $s(g)$ is regarded in (\ref{eq:G-explicit-formula}) as a subset 
of $(\Z_2)^N$. 
\end{lemma}

\begin{props} \cite[\S 1--2]{MR735406}
{\it Let $E$ be a CM field of degree $2N = [E:\Q]$. 
The group $G:= {\rm Gal}(E^{\rm nc}/\Q)$ acts on the left transitively 
on the set of $2N$ embeddings of $E$, ${\rm Hom}_{\rm field}(E,\C)$. 
When one labels the $2N$ embeddings 
by $\{ \phi_1, \bar{\phi}_1, \cdots, \phi_N, \bar{\phi}_N\}$ in a way 
the complex conjugation pairs of embeddings are respected in this labeling, 
the Galois group $G$ is regarded as the subgroup of ${\rm Im}(N,2)$ for 
some appropriate triple $(G_0, (\Z_2)^v, s)$.}
\end{props}
We call the triple a {\it Dodson data} of the CM field $E$. 
Let us call a way to label the $2N$ embeddings as above a 
{\it layout of the embeddings}. The Dodson data of $E$ depends 
on the choice of a layout of the embeddings. 
\begin{notn}
Let $(V, h)$ be a pure weight-$m$ rational Hodge structure, $V$ 
its underlying vector space over $\Q$, and $E$ a CM field within 
the algebra of Hodge endomorphisms of $(V, h)$ such that $\dim_\Q V = [E:\Q]$. 
The $2N:=[E:\Q]$ embeddings of the field $E$ in $\C$
are divided into the mutually exclusive subsets $\amalg_{p=0}^{m} \Phi^{(p,m-p)}$,
based on which eigenstates of the action of $E$ is in which Hodge components.
This structure is carried over to the set of $2N$ elements 
$\{ \phi_1, \cdots \bar{\phi}_N\} \cong G/H$, where 
$H := {\rm Gal}(E^{\rm nc}/\phi_1(E))$.  
The subgroup of ${\rm Im}(N,2)$ of permutations that map elements 
only within each one of the mutually exclusive subsets is denoted 
by $\widetilde{S}$ in general. 

When the Hodge structure $(V, h)$ is of abelian variety type 
($m=1$, $h^{1,0}=h^{0,1}=N$), the subgroup $\widetilde{S}$ is denoted 
by $\widetilde{S}_{\rm abl}$. For a K3-type [resp. CY3-type] Hodge structure 
with $m=2$, $h^{2,0}=h^{0,2}=1$ and $h^{1,1}=2N-2$ [resp. 
$m=3$, $h^{3,0}=h^{0,3}=1$ and $h^{2,1}=h^{1,2}=(N-1)$], the subgroup
is denoted by $\widetilde{S}_{\rm K3}$ [resp. $\widetilde{S}_{\rm CY3}$]. 
\end{notn}

With all the preparations, it is now clear that we should be interested 
in classifying the possible choices of Dodson data $(G_0, (\Z_2)^v, s)$ 
for $G \subset {\rm Im}(N,2)$ modulo conjugation (the adjoint action)
by $\widetilde{S}$ of appropriate type.\footnote{
Although one can work in the 
abstract group-theory language for this classification problem, it 
remains to be a separate question 
(e.g., \cite[p.7 and \S5]{MR735406}) whether there exists a CM field $E$
for a given triple $(G_0, (\Z_2)^v, s)$. See \cite{MR735406} and references 
therein for more information.
}

For a simple rational Hodge structure $(V,h)$ of CM-type, one can 
always find a CM field $E$ such that $[E:\Q] = \dim_\Q V$. The 
level-$n$ subspace of a weak-CM abelian variety always has a 
simple rational Hodge structure (by definition). So we can apply 
the thinking framework explained so far without a little modification. 
\begin{exmpl}
  \label{exmpl:K3-classify-Dodson}
For $N=2$, there are just three $\widetilde{S}_{\rm K3}$-conjugacy classes 
of Dodson data. The cases (A, B, C) at the beginning of 
section \ref{sec:abel} 
(from \cite[pp.\ 64--65, Ex.\ 8.4.(2)]{shimura2016abelian}) are 
the CM fields that have those three classes of Dodson data.\footnote{
Dodson data of a CM field $E$ with $2N = [E:\Q]$ may also be 
classified modulo conjugation by $\widetilde{S}_{\rm abl}$. Once again 
we have three distinct classes, represented by the cases (A, B, C) in \cite[Ex.\ 8.4.(2)]{shimura2016abelian}.
} %
\end{exmpl}
\begin{exmpl}
  \label{exmpl:CY3-classify-Dodson}
For $N=3$, there are eight distinct $\widetilde{S}_{CY3}$-conjugacy classes 
of Dodson data.\footnote{
See \cite[Thm. (p.7)]{MR735406} for a classification under an equivalence 
relation that differs from the conjugation by $\widetilde{S}_{\rm CY3}$ 
or $\widetilde{S}_{\rm abl}$.
The classification by $\widetilde{S}_{\rm abl}$ is in Ex.\ \ref{exmpl:abl-classify-Dodson}.
}
Three are represented by the triples of the form
$(\Z_3, 1, s)$,
with the trivial cocycle $s$ and two non-trivial\footnote{
The two non-trivial cocycles $s$ referred to here are still trivial 
in $H^1(G_0; (\Z_2)^N/(\Z_2)^v)$. 
} %
 cocycles $s$ that are not equivalent under $\widetilde{S}_{CY3}$; 
three more are represented by triples of the form
$(S_3, 1, s)$, with three different kinds of cocycles $s$ as in the
previous three. The two other $\widetilde{S}_{\rm CY3}$-conjugacy classes
are represented by the triples $(\Z_3, 3, {\rm triv}.)$ and
$(S_3, 3, {\rm triv}.)$, where ${\rm triv}.$ means $s=0$.\\
{\bf proof:} by explicit computations. \qed
\end{exmpl}

When it comes to the classification of CM-type Hodge structures $(V,h)$
that is not necessarily simple, we just have to repeat such study
for $2N \leq \dim_\Q V$ and exhaust all the partitions
$2\sum_i N_i = \dim_\Q V$.

\begin{exmpl}
\label{exmpl:abl-classify-Dodson}
Consider classification of CM-type Hodge structures of abelian variety type by using the Dodson data.
We focus on the case $\dim_\Q V=2n=6$.
Objects to be classified are (i) and (ii) combined, 
where (i) is CM-type rational Hodge structures of $V=H^1(A;\Q)$ of abelian threefolds $A$, and (ii) is choices of a commutative subalgebra 
${\cal K}$ in the Hodge endomorphism algebra of $H^1(A;\Q)$ where ${\cal K}$ 
is of $2n$ dimensions over $\Q$ and is isomorphic to a direct sum of 
CM fields, $\oplus_i K_i$. 
\begin{itemize}
\item ${\cal K}=K$ is a CM-field of degree-$(2N=6)$: this possibility is
  exhausted by classifying the $N=3$ Dodson data by the conjugation under
  $\widetilde{S}_{\rm abl}$. There are six conjugacy classes then (by explicit
  computations).  The conjugation group $\widetilde{S}_{\rm abl}$ is larger
  than $\widetilde{S}_{\rm CY3}$, and as a result the two triples of the form
  $(\Z_3,1,{\rm non.triv.})$ in Ex.\ \ref{exmpl:CY3-classify-Dodson} are conjugate to each other under
  $\widetilde{S}_{\rm abl}$; the same is true for the two triples of the form
  $(S_3,1,{\rm non.triv.})$.
\item ${\cal K} = K_2 \oplus \Q(\sqrt{p})$, where $K_2$ is a degree-$(2N=4)$
  CM field and $\Q(\sqrt{p})$ an imaginary quadratic field: this possibility 
  is further classified into four.
  Two of them are when the degree-4 CM field $K_2$ is in the case either (B)
  or (C). The other two (A.iso) and (A.noniso) are when $K_2$ is in the
  case (A), and its reflex field is either isomorphic to $\Q(\sqrt{p})$ or not. 
\item ${\cal K} = \Q(\sqrt{p_1}) \oplus \Q(\sqrt{p_2}) \oplus \Q(\sqrt{p_3})$:
  this possibility is further classified into three. (i) the three imaginary quadratic fields are all isomorphic, (ii) two of the imaginary quadratic fields are mutually isomorphic and the other one is not, and (iii) none of the three imaginary quadratic fields are isomorphic to each other. 
\end{itemize}
Even when an abelian variety $A$ under consideration is isogenous to $E^3$
of a CM elliptic curve $E$, the Dodson data of a choice of the algebra ${\cal K}$ can
be in either one of the conjugacy classes of
$(\Z_3, 1, {\rm triv}.)$, $(S_3,1,{\rm triv}.)$, (A.iso) and (i).
When $A$ is isogenous to $E_1^2 \times E_2$ where $E_1$ and $E_2$
are mutually non-isogenous CM elliptic curves, the Dodson data of a chosen algebra
${\cal K}$ may be in the conjugacy class of either (A.noniso) or (ii). 
\end{exmpl}

The classification above, which ended up with the 6+4+3=13 different kinds 
of (conjugacy classes of) Dodson data, was for the weight-1 rational Hodge 
structure of CM-type abelian threefold $A$ (and an algebra ${\cal K}$ 
acting on $H^1(A;\Q)$). For each of the 13 cases, we can compute 
the corresponding $\widetilde{S}_{\rm CY3}$-conjugacy class of Dodson data
for the rational Hodge structure on the level-$n$ subspace
$[H^n(A;\Q)]_{\ell = n}$; one just has to 
combine \cite[\S8]{shimura2016abelian} with \cite[\S1--2]{MR735406}
and generalize a little bit to do so. 
\begin{itemize}
\item $\widetilde{S}_{\rm abl}(\Z_3,1,{\rm triv.})$ and 
$\widetilde{S}_{\rm abl}(S_3,1,{\rm triv}.)$ on the weight-1 $H^1(A;\Q)$: 
$n'=1$ on the level-$n$ subspace then. 
$\widetilde{S}_{\rm abl}(\Z_3,1,{\rm non.triv.})$ on weight-1:
$\widetilde{S}_{\rm CY3}(\Z_3,1,d\epsilon_1)$ on the level-3 subspace then
 ($d\epsilon_1$ is a non-trivial cocycle whose details 
we do not write here).
$\widetilde{S}_{\rm abl}(S_3,1,{\rm non.triv.})$ on weight-1:
$\widetilde{S}_{\rm CY3}(S_3,1,d\epsilon_1)$ on the level-3 subspace then.
$\widetilde{S}_{\rm abl}(\Z_3,3,{\rm triv.})$ on weight-1: then $n'=4$ on 
the level-3 subspace, with a representative Dodson data in the form 
$(G_0,v,s) = (A_4, 1, {\rm non.triv.})$. 
$\widetilde{S}_{\rm abl}(S_3,3,{\rm triv.})$ on weight-1: then $n'=4$ on 
the level-3 subspace, with a representative Dodson data in the form 
$(G_0,v,s) = (S_4,1,{\rm non.triv.})$.  \\
The upper bound \cite[Prop.\ 1.9.1]{MR257031}, $2n' \leq 2^n$, 
is observed in these examples. 
\item (B, C) on the weight-1 $H^1(A;\Q)$: $\phi'_{(3,0)}(K')$ is then 
the composite of $(K_2)^r$ and $\Q(\sqrt{p})$ on the level-3 subspace and 
$n'=4$; $\phi'_{(3,0)}(K')$ is Galois in the case (B) and non-Galois in the 
case (C). 
(A.noniso) on weight-1: then $n'=2$ on the level-3 subspace, with 
the CM field in the case (A) in Ex.\ \ref{exmpl:K3-classify-Dodson}.
(A.iso) on weight-1: then $n'=1$ on the level-3 subspace. 
\item (i) on weight-1: then $n'=1$ on the level-3 subspace. 
(ii) on weight-1: then $n'=2$ on the level-3 subspace with $K'$ in the 
case (A) in Ex.\ \ref{exmpl:K3-classify-Dodson}.
(iii) on weight-1: then 
$K' = \Q(\sqrt{p_1}, \sqrt{p_2}, \sqrt{p_3})$ on the level-3 subspace 
(so $n'=4$).
\end{itemize}
This leads us to the following result:
\begin{props}
Let $A$ be a weak CM-type abelian variety of complex dimension $n=3$.
The CM field $K'$ and its action on the level-$n$ subspace
$[H^n(A;\Q)]_{\ell =n}$ is classified by the reflex degree $2n'$, first,
and then by the $\widetilde{S}_{\rm CY3}$-conjugation of the Dodson data
of $K'$ action on the simple rational Hodge structure on the level-$n$
subspace. There are one conjugacy class with the reflex degree $(2n'=2)$,
three conjugacy classes with $(2n'=4)$ (Ex.\ \ref{exmpl:K3-classify-Dodson}), and
eight classes with $(2n'=6)$ (Ex.\ \ref{exmpl:CY3-classify-Dodson}). 
The reflex degree can also be $2n'=8,10, 12, 14$ in principle,
because $\dim_\Q(H^3_\mathrm{prim}(T^6;\Q)) = 20-6$. 

{\it If $A$ is of CM-type (i.e., the rational Hodge structure on
  $H^1(A;\Q)$ is of CM-type), then the $\widetilde{S}_{\rm CY3}$-conjugacy
  class of the reflex field $K'$ of the CM pair of $H^1(A;\Q)$
  is one of the followings: the $(2n'=2)$ case, the $(2n'=4)$ case with $K'$ 
  in the case {\rm (A)}, the $(2n'=6)$ cases with the Dodson data in 
   $\widetilde{S}_{\rm CY3}(\Z_3,1,{\rm non.triv.})$ or
   $\widetilde{S}_{\rm CY3}(S_3,1,{\rm non.triv.})$, and   
  the $(2n'=8)$ cases with five $\widetilde{S}_{\rm CY3}$-conjugacy 
 classes\footnote{
 These $(2n'=8)$ cases are from the weight-1 Dodson data (iii), (B), (C), $\widetilde{S}_\mathrm{abl}(\Z_3,3,\mathrm{triv.})$, and $\widetilde{S}_\mathrm{abl}(S_3,3,\mathrm{triv.})$.
 The Galois groups for the level-3 subspaces of the five cases are $\Z_2\times\Z_2\times\Z_2$, $\Z_2\times\Z_4$, $\Z_2\times(\Z_4\rtimes\Z_2)$, $\Z_2\rtimes A_4$, and $\Z_2\rtimes S_4$, respectively.
 } of the Dodson data. 
}
\end{props}

\begin{rmk}
If there were a weak CM-type abelian variety $A$ of complex dimension $(n=3)$ characterized by
the reflex degree $(2n'=4)$ and $K'$ in the case (B, C), that would have 
been an example of a weak CM-type variety that is not of strong CM-type. Thm.\ \ref{thm:intro-explct} 
(or Prop.\ \ref{props:caseBC-level2-HdgStr}) indicates
that there is in fact no such abelian variety, because $n$ is odd. 

It remains to be a question if there is a weak CM-type abelian threefold $A$ with the reflex degree $(2n'=6)$ and the Dodson data in 
one of the $(8-2=6)$ $\widetilde{S}_{\rm CY3}$-conjugacy classes that cannot 
be realized by a strong CM-type abelian threefold $A$, or with 
the reflex degree greater than 8.  
\end{rmk}
    
%%%%%%%%%%%%%%%%%%%%%%%%%%%%%%%%%%%%%%%%
\subsection{A Few Calabi--Yau Threefolds of Weak CM-type}
\label{ssec:exmpl-CY}
%%%%%%%%%%%%%%%%%%%%%%%%%%%%%%%%%%%%%%%%

Besides the abelian threefolds (in section \ref{ssec:open-abel})
and ${\rm K3} \times T^2$ in section \ref{sec:K3xT2}, a Calabi--Yau 
threefold $X$ is also a projective $(n=3)$-dimensional variety with 
a trivial canonical bundle. With $h^{1,0}(X)=h^{2,0}(X)=0$ for a 
Calabi--Yau threefold, $X$ being of CM-type means that $X$ is of strong CM-type. 
We may still ask whether $X$ being of weak CM-type implies that $X$ is 
of (strong) CM-type; we do not have a logic to infer that the answer is yes;
just a few well-known examples of weak CM-type Calabi--Yau threefolds are collected here. 

\begin{exmpl}
The Fermat quintic Calabi--Yau threefold $X$ is a quintic hypersurface 
of $\C P^4$ given by the Fermat type equation $\sum_{i=1}^5 X_i^5=0$, 
where $X_{1,2,3,4,5}$ are the homogeneous coordinates of $\C P^4$. 
It is known that the level-3 subspace $[H^3(X;\Q)]_{\ell=3}$ is of 
4-dimensions within the $(b_3(X)=204)$-dimensional $H^3(X;\Q)$. In this 
example, $X$ is not only weak CM (the rational Hodge structure on 
$[H^3(X;\Q)]_{\ell=3}$ is CM with $K' \cong \Q(e^{2\pi i/5})$), but the 
remaining $[H^3(X;\Q)]_{\ell < 3}$ is also CM. 
$X$ is of strong CM-type. See \cite[Thms.\ I and II]{shioda1979hodge} (and also \cite{shioda1982geometry}).

The cyclic covering constructions in \cite{MR2129008} and \cite{rohde2009cyclic}
% Rohde = cyclic covering
% MR2129008 = Viehweg Zuo
are also designed to provide examples of Calabi--Yau threefolds of 
strong CM-type.
%    MR2905137 = Bini and van Geemen
\end{exmpl}

\begin{exmpl}A Borcea--Voisin orbifold $X$ in its simplest form is 
to use a K3 surface $S$ with a holomorphic non-symplectic 
involution $\imath_S$ and an elliptic curve $E$; using the 
holomorphic non-symplectic involution $\imath_E$ on $E$, $X$ is 
given as the canonical resolution of the $\C^2/\Z_2$ singularity 
of the $\Z_2$ quotient of $(S \times E)$ by $(\imath_S, \imath_E)$. 
The rational Hodge structure on $H^3(X;\Q)$ contains a substructure 
supported on (the $\Z_2$-invariant) $T_S\otimes H^1(E;\Q)$, where 
$T_S$ is the transcendental lattice of $S$ 
\cite{AST_1993__218__273_0}, \cite{MR1416355}. 
%
% AST_1993__218__273_0 =  Voisin mirror 
% MR1416355 = Borcea mirror construction
%
The rational Hodge substructure 
on this subspace is of CM-type when both $S$ and $E$ are of CM-type, and this subspace contains the 
level-3 subspace \cite{borcea1998calabi}.
So, such a Borcea--Voisin orbifold $X$ is 
of weak CM-type.

The orthogonal complement of the subspace $T_S \otimes H^1(E;\Q)$ 
is empty when the set of $\imath_S$-fixed points does not contain 
a curve of genus $g \geq 1$; this means that such an $X$ 
is also of strong CM-type \cite{borcea1998calabi}.
Examples of CM K3 surfaces $S$ are known \cite{MR3228297}
(for each one of Nikulin's list of families of K3 surfaces with 
an involution), even in families with 
a curve $\Sigma$ (or curves) of $\imath_S$-fixed points with 
$g(\Sigma)\geq 1$, where the rational Hodge structure on $H^1(\Sigma;\Q)$
(and hence that on the orthogonal complement in $H^3(X;\Q)$) is also 
of CM-type (and hence $X$ is strong CM).
% MR3228297 = Goto Livne Yui

It seems to be an open question in this construction (where $S$ is 
assumed to be CM) at this moment whether $X$ can be non-CM, or equivalently, 
whether the rational Hodge structure on $H^1(\Sigma;\Q)$ can be non-CM. Nothing prevents 
such a CM-type $S$ with an involution from having such a fixed point 
curve $\Sigma$, but it would be nice if we had an explicit indication and 
example where $\Sigma$ is not CM, and $X$ is not CM, or a proof that $H^1(\Sigma;\Q)$ is always CM.
\end{exmpl}

%%%%%%%%%%%%%%%%%%%%%%%%%%%%%%%%%%%%%%%%%%%%%%%%%%%%%%%%%%%%%%%%%%%%%%%%%%%%%
% \begin{figure}[tbp]
% \begin{center}
%    \includegraphics[width=.8\linewidth]{hyperelliptic_2}  
%  \caption{\label{fig:xxxx} 
%  }
%  \end{center}
% \end{figure}
%%%%%%%%%%%%%%%%%%%%%%%%%%%%%%%%%%%%%%%%%%%%%%%%%%%%%%%%%%%%%%%%%%%%%%%%%%%%%%

 \appendix 

%%%%%%%%%%%%%%%%%%%%%%%%%%%%%%%%%%%%%%%%%%%%%%%%%%%%%
\section{Various Definitions for CM-type}
 \label{sec:xxxx}
%%%%%%%%%%%%%%%%%%%%%%%%%%%%%%%%%%%%%%%%%%%%%%%%%%%%%

Multiple independent definitions have been given to the notion {\it CM-type}. 
A brief note in this appendix explains those definitions along with some 
examples, and quotes logical implications among them known in the literature.
This is a summary note than a lecture note or a review article;
readers looking for proofs of the statements here or for more systematic
expositions are referred to the literatures quoted at the beginning 
of \ref{ssec:CM-EndAlg} and \ref{ssec:CM-MTHdg}. 
This appendix does not add anything new. 
The main text is readable without this appendix, but some readers might 
be benefited from this note.  

A quick summary is that it is enough to know only the definition of CM-type rational 
Hodge structures in Def.\ \ref{def:End-CM-Hdg-str}, and that this definition 
is often the easiest to deal with in practical computations. 
Alternative (but equivalent) definitions include 
Defs.\ \ref{def:MTHdg-CM} (general) 
and \ref{def:conv-CM-abelianV} (abelian varieties) 
and Rmk.\ \ref{rmk:wgt3-equiv-cond-CM-HdgStr} (weight-3 Hodge structures).
Scattered remarks (Rmks.\ \ref{rmk:geomCM-vs-rHScm}, 
\ref{rmk:ellC}--\ref{rmk:CY3}, \ref{rmk:varHS} 
(geometry vs Hodge structure), and 
Rmk.\ \ref{rmk:equiv2-2Defs-EndAlg-MTHdg} (End algebra vs Hodge group))
explain relations among the definitions in the appendix and in the main 
text in a language a little more intuitive than elsewhere. 

%%%%%%%%%%%%%%%%%%%%%%%%%%%%%%%%%%%%%%%%%%%%%%%%%%%%%%%%%%
\subsection{Geometric Definition of CM-type Abelian Varieties}
\label{ssec:CM-EndAlg}
%%%%%%%%%%%%%%%%%%%%%%%%%%%%%%%%%%%%%%%%%%%%%%%%%%%%%%%%%%

References for the materials in the 
appendices \ref{ssec:CM-EndAlg}--\ref{ssec:CM-EndAlg2} include 
\cite[\S 1, 3, 5, 6, 8]{shimura2016abelian}, \cite[\S3]{MR3586372}, and \cite[\S 1--3]{milnecm}. 

We begin with the most traditional definition of the notion of CM-type 
for abelian varieties by using geometry. Let $A$ be a complex abelian 
variety of $n$-dimensions, with its complex analytic presentation 
$A \cong \C^n/\Lambda$, where $\Lambda$ is a free abelian group 
of rank-$2n$ embedded in $\C^n$. Then we may define two algebraic objects 
for $A$:
\begin{align}
 {\rm End}(A) & \; := \left\{ {\rm hol.~map} \; 
    \varphi: A\longrightarrow A \; | \; {\rm also~a~group~homomorphism~of~}A
   \right\}, \nonumber \\
  & \; = \left\{ \varphi \in M_{n\times n}(\C) \; | \; \varphi(\Lambda) \subset \Lambda \right\}, \\
 {\rm End}_\Q(A) & \; := \left\{ \varphi \in M_{n\times n}(\C) \; | \; 
     \varphi(\Lambda \otimes_\Z \Q) \subset \Lambda \otimes_\Z \Q \right\}. 
\end{align}
${\rm End}(A)$ is a ring, while ${\rm End}_\Q(A)$ is an algebra 
over $\Q$. 
\begin{defn}
\label{def:conv-CM-abelianV}
A complex abelian variety $A$ of $n$-dimensions is said to be of {\it CM-type}
when ${\rm End}_\Q(A)$ contains a commutative algebra ${\cal K}$ over $\Q$ 
of dimension $[{\cal K}:\Q]$ equal to $2n$. 
\end{defn}
\begin{exmpl}
\label{ex:ExE-inApp}
Think of an abelian surface $A = E \times E$ where $E$ is an elliptic curve
$\C/(\Z + i \sqrt{|p|}\Z)$ for some negative rational number $p$. 
In this example, ${\rm End}_\Q(A) = M_{2\times 2}(\Q(i\sqrt{|p|}))$. This 
non-commutative algebra contains such commutative subalgebras as 
\begin{align*}
 {\cal K} = \diag \left(\Q(i\sqrt{|p|}), \Q(i\sqrt{|p|}) \right) 
\end{align*}
and 
\begin{align*}
 {\cal K} = {\rm Span}_\Q \left\{ {\bf 1}_{2\times 2}, \; 
  i\sqrt{|p|} {\bf 1}_{2\times 2}, \;
   \left( \begin{array}{cc} & x \\ 1 & \end{array} \right), \;
   i\sqrt{|p|} \left( \begin{array}{cc} & x \\ 1 & \end{array} \right) \right\}
\end{align*}
for any $x \in \Q_{\neq 0}$. All of those commutative subalgebras 
are of dimension $4$ over $\Q$. So, this abelian surface is of CM-type
(by Def.\ \ref{def:conv-CM-abelianV}). 
\end{exmpl}
\begin{rmk}
The following fact is known: in the algebra ${\rm End}_\Q(A)$ of an   
abelian variety $A$ of complex dimension $n$, any commutative subalgebra 
${\cal K}$ over $\Q$ cannot have a dimension larger than $2n$. So, one may replace 
the condition $[{\cal K}:\Q]=2n$ in Def.\ \ref{def:conv-CM-abelianV}
by $[{\cal K}:\Q] \geq 2n$ without changing the Definition. 

\end{rmk}
%
%%%%%%%%%%%%%%%%%%%%%%%%%%%%%%%%%%%%%%%%%%%%%%%%%%%%%%%%%
\subsection{CM-type Hodge Structures Defined by the 
Endomorphism Algebras}
\label{ssec:CM-EndAlg2}
%%%%%%%%%%%%%%%%%%%%%%%%%%%%%%%%%%%%%%%%%%%%%%%%%%%%%%%%%

Note, for an abelian variety $A$, that 
${\rm End}_\Q(A) \cong {\rm End}_{\rm Hdg}(H_1(A;\Q))$; the condition 
$\varphi(\Lambda \otimes \Q) \subset \Lambda \otimes \Q$ for 
$\varphi \in {\rm End}_\Q(A)$ implies that 
$\varphi \in {\rm End}(H_1(A;\Q))$, and the condition that 
$\varphi: A \rightarrow A$ is holomorphic implies that 
the holomorphic [resp.\ anti-holomorphic] tangent vectors in $H_1(A;\C)$ 
are mapped by $\varphi$ to holomorphic [resp.\ anti-holomorphic] 
tangent vectors. The pull-back by $\varphi$'s constitute the algebra 
${\rm End}_{\rm Hdg}(H^1(A;\Q))$, where the subscript ${\rm Hdg}$ implies 
that the (1,0)-forms [resp.\ (0,1) forms] in $H^1(A;\C)$ are mapped 
by $\varphi^*$ to (1,0)-forms [resp.\ (0,1)-forms. 
So, Def.\ \ref{def:conv-CM-abelianV} is equivalent to the following 
\begin{defn}
\label{def:EndA-forH1-CM-abelianV}
A complex abelian variety $A$ of $n$-dimensions is said to be of 
{\it CM-type} when the algebra ${\rm End}_{\rm Hdg}(H^1(A;\Q))$ 
contains a commutative subalgebra ${\cal K}$ over $\Q$ with 
$[{\cal K}:\Q] = 2n$.  
\end{defn}

Def.\ \ref{def:conv-CM-abelianV} involves the algebra ${\rm End}_\Q(A)$ 
that is defined specifically for abelian varieties, while 
Def.\ \ref{def:EndA-forH1-CM-abelianV} involves only a polarizable 
rational Hodge structure on a vector space $H^1(A;\Q)$. 
As a preparation for applications to varieties other than 
abelian varieties, more jargons are explained in 
Defs.\ \ref{def:rHdgStr}--\ref{def:pol-rat-Hstr}. 
\begin{defn}
\label{def:rHdgStr}
A pair $(V_\Q, \phi)$ of a finite dimensional vector space $V_\Q$ 
over $\Q$ and an isomorphism (as vector spaces over $\C$)
\begin{align}
  \phi: V_\Q \otimes \C \cong \oplus_{p,q \in \Z}^{(p+q=m)} V^{p,q}
\end{align}
is called a {\it rational Hodge structure of weight-$m$} when 
% the following conditions are satisfied. First, 
the complex conjugation 
operation acting on the $\C$ tensor factor of $V_\Q \otimes \C$
converts $V^{p,q}$ to $V^{q,p}$. 
We may allow ourselves sometimes to use an expression ``a rational 
Hodge structure $V_\Q$'' when the choice of $\phi$ in $(V_\Q, \phi)$ 
is obvious from contexts. 

A {\it Hodge morphism} $f$ from a rational Hodge structure $(V_\Q, \phi)$
to $(V'_\Q, \phi')$ of the same weight is a vector space homomorphism 
$f: V_\Q \rightarrow V'_\Q$ that maps $V^{p,q}$ to $(V')^{p,q}$. 
\end{defn}
\begin{defn}
A rational Hodge structure $(V_\Q, \phi)$ is said to have a 
{\it rational Hodge substructure} $(W_\Q, \phi|_W)$ when there is 
a vector subspace $W_\Q$ over $\Q$ within $V_\Q$ such that 
\begin{align}
  W^{p,q} := \phi (W_\Q \otimes \C)\cap V^{p,q}, \qquad 
  \oplus_{p,q}^{(p+q=m)} W^{p,q} = \phi(W_\Q \otimes \C). 
\end{align}
A rational Hodge structure $(V_\Q,\phi)$ is said to be {\it simple}, if it has no rational Hodge substructure besides the zero vector space and itself.

\end{defn}
\begin{defn}
\label{def:pol-rat-Hstr}
A triple $(V_\Q, \phi; {\cal Q})$ is called a {\it polarized rational Hodge 
structure} when (i) $(V_\Q, \phi)$ is a rational Hodge structure of some 
weight-$m$, (ii) ${\cal Q} : V_\Q \times V_\Q \rightarrow 
\Q$ a bilinear form that is symmetric [resp.\ anti-symmetric] when $m$ is even 
[resp.\ odd], and (iii) 
the following two conditions are satisfied:
\begin{align}
  \chi \in V^{p,q}, \; \psi \in V^{p',q'} \quad & \; {\rm then~} 
   {\cal Q}(\chi,\psi) \neq 0 \quad {\rm only~if~} p+p'=q+q'=m, 
    \label{eq:cond-polarization-w-J} \\
  \chi \in V^{p,q} \backslash \{0\} \quad & \; {\rm then~}
    (-1)^{\frac{m(m+1)}{2}}{\cal Q}(\chi, J^* \bar{\chi}) > 0. 
    \label{eq:cond-polarization-pDef}
\end{align}
Here, $J^*$ is the operator multiplying $i^{p-q}$ to $\psi \in V^{p,q}$. 
The last condition is equivalent to the positive definiteness of 
the symmetric bilinear form $H_{\cal Q}(-,-):= (-1)^{m(m+1)/2}{\cal Q}(-,J^*-)$
on $V_\Q \otimes \R$.
A bilinear form ${\cal Q}$ satisfying the conditions (ii) and (iii) is called 
a polarization of $(V_\Q, \phi)$. 

A rational Hodge structure $(V_\Q, \phi)$ is said to be {\it polarizable}, 
when there exists a polarization ${\cal Q}$ such that $(V_\Q, \phi; {\cal Q})$ 
becomes a polarized rational Hodge structure. 

A bilinear form ${\cal Q}$ satisfying the condition (ii), 
(iii)-(\ref{eq:cond-polarization-w-J}) is called a {\it polarization 
of index} $k$, if the symmetric bilinear form $H_{\cal Q}$ is non-degenerate, 
and has $k$ negative eigenvalues; polarizations defined already are 
polarizations with index 0 (following \cite{birkenhaketori}). 
\end{defn}
\begin{rmk}
\label{rmk:HdgStr-4-projVar}
Let $X$ be a complex $n$-dimensional projective non-singular variety. 
Then the vector 
space $H^k(X,\Q)$ (with $k=0,1,\cdots, 2n$) is given a rational Hodge 
structure that is determined uniquely by the complex structure of $X$. 
The rational Hodge structure on $H^k(X;\Q)$ is known to be polarizable. 
An embedding of $X$ into a projective space determines a choice of 
${\cal Q}$.  
%
% (iii) there exists rational Hodge substructures $(V_\Q^{(s)}, \phi^{(s)})$
% of $(V_\Q, \phi)$ labeled by $s \in \Z_{\geq 0}$ such that 
% %
% \begin{align}
%   {\cal Q}(\psi, \chi) = 0 \qquad {\rm when~}
%     \psi \in V^{(s)}, \quad \chi \in V^{(t)}, \quad s \neq t, 
% \end{align}
% %
% (iv) $V_\Q \cong \oplus_s V_\Q^{(s)}$, 
% 
\end{rmk}
\begin{rmk}
\label{rmk:forAbelV-equiv-geomDef-H1end}
For an abelian variety $A$, the algebra ${\rm End}_{\rm Hdg}(H^1(A;\Q))$ 
is the algebra of Hodge morphisms of the rational Hodge structure on 
$H^1(A;\Q)$. So, we may read the 
Defs.\ \ref{def:conv-CM-abelianV}=\ref{def:EndA-forH1-CM-abelianV} as 
imposing conditions on the rational Hodge structure. 
\end{rmk}
As a generalization of Def.\ \ref{def:EndA-forH1-CM-abelianV}, 
we may introduce the following notions:
\begin{defn}
\label{def:End-CM-Hdg-str}
Let $(V_\Q, \phi)$ be a rational Hodge structure. It is said to have 
{\it sufficiently many complex multiplications}, when the algebra 
of Hodge morphisms ${\rm End}(V_\Q, \phi)$
(an equivalent notation ${\rm End}_{\rm Hdg}(V_\Q)$)
contains a commutative semi-simple subalgebra ${\cal K}$ over $\Q$ with 
$[{\cal K}:\Q] = \dim_\Q V_\Q$. 

Let $(V_\Q, \phi)$ be a rational Hodge structure that is polarizable. 
It is said to be of {\it CM-type}, when it has sufficiently many 
complex multiplication. 
\end{defn}
\begin{rmk}
Note that the algebra ${\rm End}(V_\Q, \phi)$ depends only on a pair 
$(V_\Q, \phi)$, not on a choice of a polarization ${\cal Q}$ even when 
a rational Hodge structure $(V_\Q, \phi)$ is polarizable. It is just a
matter of conventions (adopted in many literatures) to reserve the 
term CM-type only for rational Hodge structures that are polarizable (this is for a good reason; a lot more constraints on the algebra 
${\rm End}(V_\Q,\phi)$ are available when there is a polarization).
\end{rmk}
\begin{rmk}
For a simple CM-type polarizable Hodge structure, the subalgebra ${\cal K}$ found in Def.\ \ref{def:End-CM-Hdg-str} coincides with ${\rm End}(V_\Q,\phi)$ in fact, and ${\cal K}$ is isomorphic to a CM field.
\end{rmk}
\begin{rmk}
\label{rmk:geomCM-vs-rHScm}
For a complex projective non-singular variety $X$ of $n$-dimensions, not 
necessarily an abelian variety, one may combine Def.\ \ref{def:End-CM-Hdg-str}
with Rmk.\ \ref{rmk:HdgStr-4-projVar} to introduce the notions of CM-type 
to the variety $X$ itself. There are multiple different ways to do so, 
however, because one might apply Def.\ \ref{def:End-CM-Hdg-str} to 
various Hodge substructures of $H^*(X;\Q)$. We have presented some of them as 
Def.\ \ref{defn:strong-CM}, \ref{defn:prim-CM} and \ref{defn:weak-CM}. 
Precise relation among them is the main topics of this article. 
\end{rmk}
\begin{rmk}
\label{rmk:ellC}
When $X$ is an abelian variety of $n$-dimensions, the following is known. 
When the weight-1 rational Hodge structure on $H^1(X;\Q)$ is of CM-type 
(in the sense of Def.\ \ref{def:End-CM-Hdg-str}), then the rational Hodge 
structure on $H^k(X;\Q)$ is of CM-type for each of $k=0,1,\cdots, 2n$.  
So, $X$ is of strong CM-type. Conversely, when $X$ is of strong CM-type, 
the weight-1 rational Hodge structure on $H^1(X;\Q)$ is of CM-type obviously.
So, (Rmk.\ \ref{rmk:forAbelV-equiv-geomDef-H1end}) $X$ is of strong CM-type
if and only if $X$ is of CM-type in the conventional 
geometric definition Def.\ \ref{def:conv-CM-abelianV}. 
\end{rmk}

\begin{rmk}
\label{rmk:K3}
Let $X$ be a complex projective non-singular K3 surface. The rational Hodge 
structure is non-trivial only on $H^2(X;\Q)$. The transcendental lattice 
tensored with $\Q$, denoted by $T_X\otimes \Q$, supports a rational Hodge 
substructure of the Hodge structure on $H^2(X;\Q)$
%  the rational Hodge structure on 
% $T_X\otimes \Q$ does not have a rational Hodge proper substructure 
% (by definition on the transcendental lattice). Therefore 
The K3 surface $X$ is of strong CM-type (and also of CM-type in the sense 
of Def.\ \ref{defn:prim-CM}) if and only if the rational Hodge structure 
on $T_X \otimes \Q$ is of CM-type in the sense of 
Def.\ \ref{def:End-CM-Hdg-str}.   

Although a complex projective non-singular K3 surface of strong CM-type 
has many linear maps $\varphi$ from $T_X \otimes \Q$ to itself that preserves 
the rational Hodge structure, those $\varphi$'s are not necessarily 
associated with a map from the original geometry $X$ to itself; if the 
$\varphi$'s were maps from $T_X$ to itself (than from $T_X \otimes \Q$ 
to itself), the Torelli theorem for K3 surfaces would have had things to say 
on automorphisms of $X$ associated with $\varphi$'s. For a higher dimensional 
variety $X$ that is not an abelian variety, linear maps 
in ${\rm End}_{\rm Hdg}(H^*(X;\Q))$ cannot necessarily be translated into 
geometric 
automorphisms of the original geometry $X$.
\end{rmk}
\begin{rmk}
\label{rmk:CY3}
Let $X$ be a complex projective Calabi--Yau threefold, where 
$h^{1,0}(X)=h^{2,0}(X)=0$. The rational Hodge structure 
on $H^2(X;\Q)$ consists of only the Hodge (1,1) component, so it is obviously 
of CM-type in the sense of Def.\ \ref{def:End-CM-Hdg-str}; the same 
is true with $H^4(X;\Q)$.  
So, $X$ is of strong CM-type (or of CM-type in Def.\ \ref{defn:prim-CM})
if and only if the rational Hodge structure on 
$H^3(X;\Q) = H^3_{\rm prim}(X;\Q)$ is of CM-type in the sense of 
Def.\ \ref{def:End-CM-Hdg-str}.  
\end{rmk}
%

%%%%%%%%%%%%%%%%%%%%%%%%%%%%%%%%%%%%%%%%%%%%%%%%%%%%%%%%%
\subsection{CM-type Hodge Structures Defined by the Hodge Groups
/ Mumford--Tate Groups}
\label{ssec:CM-MTHdg}
%%%%%%%%%%%%%%%%%%%%%%%%%%%%%%%%%%%%%%%%%%%%%%%%%%%%%%%%% 

Lecture notes that discuss Hodge structures in the language 
of Mumford--Tate group/Hodge groups include 
\cite{BMoonen-Lect-MTgrp-Intro}, \cite[\S1]{rohde2009cyclic}, \cite{MR2931229}, \cite{MR3290134}, 
\cite[\S4]{milnecm}. 
Books on focused on algebraic group are also available
 \cite{malle_testerman_2011} and \cite{milne_2017}. 

There is yet another way to introduce the notion of CM-type/sufficiently 
many complex multiplications on a rational Hodge structure. 
Def.\ \ref{def:End-CM-Hdg-str} we have seen so far is 
known to be equivalent to Def.\ \ref{def:MTHdg-CM} introduced below, 
as stated in Rmk.\ \ref{rmk:equiv-2Defs-EndAlg-MTHdg}. 
Before writing down Def.\ \ref{def:MTHdg-CM}, we need a little preparations: 
the additional information $\phi$ on $V_\Q$ in a rational Hodge structure 
$(V_\Q, \phi)$ is recaptured in the language of algebraic groups 
(Rmk.\ \ref{rmk:ratHS-by-DlgnTrsRepr}), as follows:
\begin{defn}
An {\it algebraic group} $G$ over a field $k$ is a variety over $k$ along with 
the multiplication law $G \times G \rightarrow G$ and the inverse law 
$G \rightarrow G$ given as regular maps of varieties over $k$. 
The field $k$ is not necessarily algebraically closed. 
See \cite{malle_testerman_2011} and \cite{milne_2017} for more information. 
\end{defn}
\begin{exmpl}
${\rm Res}_{\C/\R}(\mathbb{G}_m)$ denotes\footnote{
More generally, for a field extension $E/F$, ${\rm Res}_{E/F}(\mathbb{G}_m)$ 
stands for the algebraic group over $F$ given by 
$\{ (x_{1,\cdots,[E:F]};z') \; | \; z' {\rm Nm}_{E/F}(x_a \omega_a) =1\}$, 
where $E= {\rm Span}_F \{\omega_{a=1,\cdots, [E:F]} \}$. 
} %
 an algebraic group over $\R$ 
that is the affine algebraic variety given by $\{ (x,y,z') \; | \; 
(x^2+y^2)z' = 1\}$, with the multiplication law given by 
\[
((x_1,y_1, z'_1) , (x_2,y_2,z'_2))
 \longmapsto (x_1x_2-y_1y_2, x_1y_2+y_1x_2, z'_1z'_2), 
\]
and the inverse law 
by $(x,y,z') \longmapsto (x z', -y z',x^2+y^2)$. 
% 
% We say that such an algebraic 
% group is defined over the field $\R$ because 
% the polynomials appearing in the multiplication law and the inverse law 
% involve only the coefficients in $\R$.

Another algebraic group $U^1$ over $\R$ is obtained by imposing 
one more equation $x^2+y^2=1$ on ${\rm Res}_{\C/\R}(\mathbb{G}_m)$. 

The group ${\rm GL}(V_\Q)$ [resp.\ ${\rm GL}(V_\R)$] of a vector space 
$V_\Q$ over $\Q$ [resp.\ $V_\R$ over $\R$] is also an algebraic group 
defined over $\Q$ [resp.\ $\R$].  
\end{exmpl}
\begin{defn}
For a pair of algebraic groups $G_k$ and $G'_k$, both defined over 
a common field $k$, a {\it morphism} $f: G \rightarrow_{/k} G'$ 
{\it of algebraic groups over $k$} is a group homomorphism 
$G \rightarrow G'$ that is also a morphism defined over $k$
when we see $G$ and $G'$ as algebraic varieties over $k$; the last 
condition means in a more colloquial terms that the polynomials 
describing the map $G \rightarrow G'$ of affine subvarieties involve 
only the coefficients in the field $k$. 
\end{defn}
\begin{rmk}
\label{rmk:ratHS-by-DlgnTrsRepr}
Let $(V_\Q, \phi)$ be a rational Hodge structure. Then one can specify 
a representation $h_\phi$ of the multiplication group $\C^\times$ on 
$V_\Q\otimes \C$ given by 
\begin{align}
  h_\phi: \C^\times \ni z \longmapsto (z^p \bar{z}^q)\times \; {\rm on} \;
   V^{p,q} \subset V_\Q \otimes \C. 
\end{align}
It is possible to think of $h_\phi(z)$ for $z \in \C^\times$ as $\R$-valued 
matrix acting on $V_\Q \otimes \R$ (not on $V_\Q \otimes \C$) in fact. 
The representation $h_\phi$ may be regarded, therefore, as a morphism\footnote{
$(x,y,z') \in {\rm Res}_{\C/\R}(\mathbb{G}_m)$ with $x,y,z' \in \R$
is regarded as a complex number $z =(x+iy) \in \C^\times$. 
} % 
\begin{align}
 h_\phi : {\rm Res}_{\C/\R}(\mathbb{G}_m) & \; \longrightarrow_{/\R}
     {\rm GL}(V_\Q \otimes \R); 
\end{align}
of algebraic groups over $\R$. Conversely, from a given pair 
$(V_\Q, h)$ of a finite dimensional vector space $V_\Q$ over $\Q$ 
and a morphism $h: {\rm Res}_{\C/\R}(\mathbb{G}_m) \rightarrow_{/\R} {\rm GL}(V_\Q \otimes \R)$ over $\R$, one may extract a Hodge decomposition 
$\phi: V\otimes \C \cong \oplus_{p,q}V^{p,q}$.  So, we will use 
the notations $(V_\Q, \phi)$ and $(V_\Q, h)$ for a rational Hodge structure 
interchangeably in the rest of this appendix. 
\end{rmk}
\begin{defn}
Let $(V_\Q, h)$ be a rational Hodge structure. 
The {\it Mumford--Tate group} of $(V_\Q, h)$, denoted by $MT(V_\Q,h)$, 
is the smallest algebraic subgroup of the algebraic group ${\rm GL}(V_\Q)$ 
defined over $\Q$ that contains $h({\rm Res}_{\C/\R}(\mathbb{G}_m))$. 
%
% whose set of $\R$-points contain the image of the $\R$-points of 
% ${\rm Res}_{\C/\R}(\mathbb{G}_m)$ by $h$. 

% Here, for an algebraic group 
% $G$ defined over a field $k$, the set of $k'$-points $G(k')$ for 
% an extension field $k'/k$ is the Affine variety $G\otimes_k k'$. 
% For example, the set of $\R$-points of ${\rm Res}_{\C/\R}(\mathbb{G}_m)$ 
% is $\C^\times$. 

%
% The set of its $\R$-points, $U^1(\R)$, is the circle of unit 
% radius in $\R^2$. 

Let $(V_\Q, h)$ be a rational Hodge structure. 
The {\it Hodge group} (or {\it special Mumford--Tate group}) of $(V_\Q,h)$, 
denoted by ${\rm Hg}(V_\Q,h)$, is the smallest algebraic subgroup of the 
algebraic group ${\rm GL}(V_\Q)$ defined over $\Q$ that contains $h(U^1)$.

Some literatures define $MT(V_\Q, h)$ or ${\rm Hg}(V_\Q,h)$ only 
for a rational Hodge structure that is polarizable. In this appendix, 
those two groups are introduced for a rational Hodge structure that 
may or may not have a polarization. 
\end{defn}
\begin{rmk}
The Mumford--Tate group and the Hodge group of a rational Hodge structure 
$(V_\Q, h)$ of weight-$m$ are connected. When $(V_\Q, h)$ has a polarization 
${\cal Q}$, then $MT(V_\Q, \phi)$ and ${\rm Hg}(V_\Q, h)$ are not just 
subgroups of ${\rm GL}(V_\Q)$ and ${\rm SL}(V_\Q)$, respectively, but 
a subgroup of 
\begin{align}
  \mathbb{G}{\rm Sp}(V_\Q, {\cal Q}) \quad {\rm resp.} \quad 
   {\rm Sp}(V_\Q, {\cal Q}), \qquad {\rm when} \;  m \; {\rm is~odd}, \\
  \mathbb{G}{\rm O}(V_\Q, {\cal Q}) \quad {\rm resp.} \quad 
  {\rm SO}(V_\Q, {\cal Q}), \qquad {\rm when} \; m \; {\rm is~even}. 
\end{align}
Here, the groups $\mathbb{G}{\rm Sp}$ and $\mathbb{G}{\rm O}$ 
consist of elements $g$ such that ${\cal Q}(gx, gy) = \lambda_g {\cal Q}(x,y)$ 
for ${}^\forall x,y \in V_\Q$ with ${}^\exists \lambda_g$ independent of $x,y$.
\end{rmk}
\begin{rmk}
\label{rmk:MT-Hdg-rltn}
For a rational Hodge structure $(V_\Q, \phi)$, it is known that 
\begin{align}
  MT(V_\Q, h) = (\mathbb{G}_m {\bf 1}_{V_\Q})\cdot {\rm Hg}(V_\Q, h), 
\end{align}
where $(\mathbb{G}_m {\bf 1}_{V_\Q})$ is the algebraic subgroup of 
${\rm GL}(V_\Q)$ of scalar multiplication on $V_\Q$.  
\end{rmk}
\begin{rmk}
\label{rmk:EndAlg-MTHdg-rltn}
Let $(V_\Q, h)$ be a rational Hodge structure. It is known that 
\begin{align}
 {\rm End}(V_\Q, h) = {\rm GL}(V_\Q)^{MT(V_\Q,h)} = {\rm GL}(V_\Q)^{{\rm Hg}(V_\Q,h)}.
\end{align}
This relation allows us to determine the Mumford--Tate/Hodge group 
from the endomorphism algebra ${\rm End}(V_\Q, h)$, or also in the 
other way around. 

The spirit of Def.\ \ref{def:End-CM-Hdg-str} is that a rational Hodge 
structure [resp.\ a polarizable rational Hodge structure] is with 
sufficiently many complex multiplications [resp.\ of CM-type] 
when the algebra ${\rm End}(V_\Q,h)$ within ${\rm End}(V_\Q)$ 
is large. This means, given Rmk.\ \ref{rmk:EndAlg-MTHdg-rltn}, 
that the groups $MT(V_\Q, h)$ and ${\rm Hg}(V_\Q,h)$ are small then. 
\end{rmk}
\begin{defn}
\label{def:MTHdg-CM}
A rational Hodge structure [resp.\ a polarizable rational Hodge structure] 
$(V_\Q, h)$ is said to be with {\it sufficiently 
many complex multiplications} [resp.\ of {\it CM-type}], when 
the group $MT(V_\Q, h)$ is an algebraic torus (see below), or equivalently, 
the group ${\rm Hdg}(V_\Q, h)$ is an algebraic torus; 
\end{defn}
here, 
\begin{defn}
an algebraic group $G$ over a field $k$ of characteristic 0 is said to 
be an {\it algebraic torus}, when the algebraic group $G \otimes_k \bar{k}$, 
where $\bar{k}$ is an algebraic closure of $k$, is isomorphic to the 
product of a finite number of copies of $\mathbb{G}_m$. Here, 
$\mathbb{G}_m$ is the algebraic group given by $\left\{ (x,z') \; | \; xz'=1
\right\}$ with the multiplication law $((x_1,z'_1), (x_2, z'_2)) \longmapsto 
(x_1x_2, z'_1z'_2)$ and the inverse law $(x,z') \longmapsto (z',x)$. 
\end{defn}
\begin{rmk}
\label{rmk:equiv-2Defs-EndAlg-MTHdg}
The two definitions, Def.\ \ref{def:End-CM-Hdg-str} and \ref{def:MTHdg-CM}, 
are known to be equivalent. 

A sketch of the argument (not a proof)  
for \ref{def:MTHdg-CM}$\Rightarrow$\ref{def:End-CM-Hdg-str} is to choose 
a maximal torus subgroup in ${\rm GL}(V_\Q)$ that contains $MT(V_\Q, h)$ 
or ${\rm Hg}(V_\Q, h)$, and set ${\cal K} \subset {\rm End}(V_\Q)$ 
as the algebra over $\Q$ generated by the rational elements of the maximal 
torus. This ${\cal K}$ is a commutative subalgebra required in 
Def.\ \ref{def:End-CM-Hdg-str}. The converse 
argument \ref{def:End-CM-Hdg-str}$\Rightarrow$\ref{def:MTHdg-CM} 
proceeds by choosing a commutative semi-simple subalgebra ${\cal K} \subset 
{\rm End}(V_\Q, h)$ of dimension $\dim_\Q V_\Q$. The group of 
invertible elements of ${\cal K}$ contains the groups $MT(V_\Q, h)$ 
and ${\rm Hg}(V_\Q,h)$, from which one can understand  
the commutativity of the multiplication law of $MT(V_\Q,h)$ and 
${\rm Hg}(V_\Q,h)$. 
%
% See \cite[footnote 17 in v.2]{Kidambi:2022wvh}
% for a little more in the case without a polarization. 
%
\end{rmk}
\begin{rmk}
\label{rmk:equiv2-2Defs-EndAlg-MTHdg}
Neither of the two equivalent definitions is regarded less  
fundamental than the other; given Rmk.\ \ref{rmk:EndAlg-MTHdg-rltn}, 
each of the two algebraic objects ${\rm End}(V_\Q, h)$ and 
${\rm Hg}(V_\Q,h)$ has complete information of the other. 
In practical computations, one might use Def.\ \ref{def:End-CM-Hdg-str}
and avoid computing the Hodge/Mumford--Tate groups, or use 
Def.\ \ref{def:MTHdg-CM} and avoid finding a commutative subalgebra ${\cal K}$. 
\end{rmk}
\begin{exmpl}
In Ex.\ \ref{ex:ExE-inApp}, all the commutative subalgebras ${\cal K}$
written there contain 
\begin{align}
 MT(V_\Q, h) = {\rm Res}_{\Q(i\sqrt{|p|})/\Q}(\mathbb{G}_m) {\bf 1}_{2\times 2}
 \subset {\rm End}(H_1(E;\Q) \oplus H_1(E;\Q)). 
\end{align}
It is known that the algebraic group ${\rm Res}_{F/\Q}(\mathbb{G}_m)$
over $\Q$ is an algebraic torus for any finite extension $F/\Q$. 
So, the rational Hodge structure on $H^1(A;\Q)$ in Ex.\ \ref{ex:ExE-inApp}
is of CM-type in the sense of Def.\ \ref{def:MTHdg-CM}. 
\end{exmpl}
%

%%%%%%%%%%%%%%%%%%%%%%%%%%%%%%%%%%%%%%%%%%%%%%%%%%%%%
\subsection{Complex Multiplications on Intermediate Jacobians}
%%%%%%%%%%%%%%%%%%%%%%%%%%%%%%%%%%%%%%%%%%%%%%%%%%%%

Ref.\ \cite{borcea1998calabi} observed that the Hodge decomposition 
of a weight-3 rational Hodge structure $(V_\phi, \phi)$, 
\begin{align}
  \phi: V_\Q \otimes \C \cong V^{3,0} \oplus V^{2,1} \oplus V^{1,2} \oplus V^{0,3}, 
\end{align}
can also be captured by the two decompositions 
\begin{align}
 \phi^W: V_\Q \otimes \C & \; \cong (V^{2,1}\oplus V^{0,3}) \oplus 
   (V^{1,2}\oplus V^{3,0}) =: V^{1,0}_{\rm Weil} \oplus V^{0,1}_{\rm Weil}, \\
 \phi^G: V_\Q \otimes \C & \; \cong (V^{3,0}\oplus V^{2,1}) \oplus 
   (V^{0,3}\oplus V^{1,2}) =: V^{1,0}_{\rm Griffiths} \oplus V^{0,1}_{\rm Griffiths}.  
\end{align}

More generally, when an integral Hodge structure $(V_\Q, \phi)$ 
of an odd weight-$m$ is given, two ways to extract weight-1 
rational Hodge structures have been discussed in the literature. 
One is the Weil intermediate Jacobian and the other is 
the Griffiths intermediate Jacobian (see \cite[\S4]{birkenhaketori}). 
When we deal with only rational Hodge structures, $(V_\Q, \phi^W)$
and $(V_\Q, \phi^G)$ are those of the Weil and Griffiths intermediate 
Jacobians, respectively. 
For an odd weight $m \geq 5$, information in the two weight-1 rational 
Hodge structures cannot reproduce the original rational Hodge structure
$(V_\Q, \phi)$.

When a rational Hodge structure $(V_\Q, \phi)$ of an odd weight 
is polarizable, then 
the weight-1 rational Hodge structure $(V_\Q, \phi^W)$ is also polarizable, 
while $(V_\Q, \phi^G)$ has a polarization of some index not necessarily zero. 
\begin{rmk}
\label{rmk:wgt3-equiv-cond-CM-HdgStr}
Let $(V_\Q, \phi)$ be a rational Hodge structure of weight-3. Then the 
two conditions are equivalent \cite{borcea1998calabi}:
\begin{itemize}
\item the rational Hodge structure $(V_\Q, \phi)$ has sufficiently many complex 
multiplications (in the sense of the equivalent Defs.\ \ref{def:End-CM-Hdg-str} 
and \ref{def:MTHdg-CM}), and 
\item both $(V_\Q, \phi^W)$ and $(V_\Q, \phi^G)$ have sufficiently many 
complex multiplications, and their Hodge groups ${\rm Hg}(V_\Q, \phi^W)$ 
and ${\rm Hg}(V_\Q,\phi^G)$ within ${\rm GL}(V_\Q)$ commute with each other. 
\end{itemize}
The condition that the two Hodge groups commute is necessary in order to 
make sure that there exists a maximal torus containing both of the Hodge 
groups. One may also use this, 
\begin{itemize}
\item there exists a commutative algebra ${\cal K}$ with $[{\cal K}:\Q]=
\dim_\Q V_\Q$ within ${\rm End}(V_\Q, \phi^W) \cap {\rm End}(V_\Q,\phi^G)$, 
\end{itemize}
as another condition equivalent to the two above.
\end{rmk}
\begin{rmk}
The previous remark can be applied to any rational Hodge substructures 
of $H^3(X;\Q)$ of a non-singular complex projective Calabi--Yau threefold 
$X$. When we apply the equivalence above only to the level-3 substructure 
$[H^3(X;\Q)]_{\ell=3}$, then $X$ being of weak CM-type means that 
both the Griffiths and Weil intermediate Jacobians of the level-3 
substructure are of sufficiently many complex multiplications and CM-type, respectively 
(and their Hodge groups commute). When $X$ is of CM-type (in the sense 
of Def.\ \ref{defn:prim-CM}), the Griffiths [resp.\ Weil] intermediate 
Jacobian of $H^3_{\rm prim}(X;\Q)$ is with sufficiently many complex 
multiplications [resp.\ of CM-type]. The previous remark only allows 
us to translate the same information on a given rational Hodge substructure
of $H^3(X;\Q)$ between different languages, but does not infer anything 
on other rational Hodge substructures of $H^*(X;\Q)$ in general. 
\end{rmk}
\begin{rmk}
\label{rmk:varHS}
Thm.\ \ref{thm:intro-main} in the main text, however, does that.
%When we apply Thm.\ \ref{thm:intro-main} to the case of an abelian threefold 
%$X$,
When $X$ is an abelian variety of weak CM-type with the reflex degree $2n'\leq 4$, all the other rational Hodge substructures of $H^*(X;\Q)$ are also 
of CM-type. 
The main theme of this article in the main text 
is how much we can learn about the various rational Hodge substructures of $H^*(X;\Q)$
from a special property on one rational Hodge substructure 
$[H^n(X;\Q)]_{\ell =n}$, when $X$ is a non-singular complex projective 
variety with a trivial canonical bundle. 

\end{rmk}
%

% When there exists a commutative subalgebra ${\cal K}$ over $\Q$ 
% within ${\rm Hom}(V_\Q)$ that has $[{\cal K}:\Q] = \dim_\Q V_\Q$, and 
% is contained both in ${\rm End}((V_\Q, \phi^{\rm Grff}))$ and 
% ${\rm End}((V_\Q, \phi^{\rm Weil}))$, then $(V_\Q, \phi)$ is also of 
% CM-type in the sense of Def.\ \ref{def:End-CM-Hdg-str}. 
% [it may be that the condition here is stronger than just demanding 
% that two rational Hodge structures $(V_\Q, \phi^{\rm Grff})$ and 
% $(V_\Q, \phi^{\rm Weil})$ of weight-1 are both of CM-type; would be nice 
% if we find an example of this subtle marginal case \cite{borcea1998calabi} 
% imposes an extra condition that $hg(Weil)$ and $hg(Grff)$ commute]

%

%  \\
%  \C^\times \ni z & \; \longmapsto (z^p\bar{z}^q \times) \quad 
%  {\rm on~}V^{p,q} \subset \oplus V^{p.q} \cong V_\Q \otimes \C.  
%

%
% \begin{defn}
% Probably an affine algebraic group that is commutative and connected 
% is a torus, and vice versa. When the condition ``connected'' is dropped, 
% then $G = {\rm Spec}(k[x]/(x^n-1)) \cong \Z/n\Z$ is affine and commutative, 
% but not a torus. When the condition ``affine'' is dropped, then 
% an abelian variety is not a torus. For more information, perhaps 
% we should examine [Milne lecture note ``algebraic group'', Ch.\ 14]
% 
% [Milne ``AG'' p.\ 72] an algebraic group that admits a faithful linear 
% representation is called a linear algebraic group; an algebraic group 
% is linear iff it is affine. 
% \end{defn}
% %

% ..............................................

%
% \begin{thebibliography}{99}

\bibliographystyle{alpha}             % alphabetical order
\bibliography{complexMultiplication}  % use complexMultiplication.bib

\end{document}